\documentclass[dvips,titlepage,A4,11pt,fleqn]{article}
\usepackage{graphicx}
\usepackage{amsmath}
\usepackage{latexsym}
\usepackage{amssymb, times}
\usepackage{float}
\usepackage{supertabular}
\usepackage{epsfig}
\usepackage[frame]{xy}
\usepackage{fancybox}
\usepackage{tabularx}
\usepackage{multirow}
\usepackage{calc}
\usepackage{color}
\usepackage{parskip}
\usepackage{endnotes}
\newcommand{\be}{\begin{equation*}}
\newcommand{\ee}{\end{equation*}}
\newcommand{\ben}{\begin{equation}}
\newcommand{\een}{\end{equation}}
\newcommand{\bea}{\begin{eqnarray*}}
\newcommand{\eea}{\end{eqnarray*}}
\newcommand{\bean}{\begin{eqnarray}}
\newcommand{\eean}{\end{eqnarray}}
\newcommand{\ba}{\begin{array}}
\newcommand{\ea}{\end{array}}
\newcommand{\lba}{\left[ \begin{array}}
\newcommand{\ear}{\end{array} \right]}

\newcommand{\cR}{{\rm I\hspace{-0.7mm}R}}
\newtheorem{theorem}{Theorem}
\newtheorem{lemma}{Lemma}
\newtheorem{corollary}{Corollary}
\newtheorem{remark}{Remark}

\newlength{\LL} 
\newlength{\NL} 
\newlength{\MD} 
\definecolor{brown}{rgb}{0.6,0.2,0.1}
\title{Derivation of a transfer function model for a high pressure pipeline}

\author{{\Large \bf Paulo Lopes dos Santos} \\ \vspace{-5mm}        \\
        Dep. de Engenharia Electrot\'{e}nica e Computadores \\
        Faculdade de Engenharia da Universidade do Porto \\
        Rua Dr Roberto Frias, s/n                         \\
        4200-464 Porto, Portugal                                     \\
       pjsantos@fe.up.pt\\[4mm]
        \Large \bf T-P de Azevedo Perdico\'{u}lis\\
        ISR-Coimbra\\
        UTAD\\
        tazevedo@utad.pt\\[4mm]
         \Large \bf J. A. Ramos\\
         Farquhar College of Arts and Sciences\\
          Division of Mathematics \\
         Science,
and Technology\\
 Nova Southeastern University\\ 3301 College Avenue\\ Fort Lauderdale, FL 33314, USA  \\jr1284@nova.edu\\[4mm]
       \Large \bf G. Jank\\
       RWTH-University of Technology\\ Department of Mathematics\\ 52056 Aachen, Germany\\
jank@math2.rwthe-aachen.de \\[4mm]
   \Large \bf J. L. Martins de Carvalho\\
       Dep. de Engenharia Electrot\'{e}nica e Computadores \\
        Faculdade de Engenharia da Universidade do Porto \\
        Rua Dr Roberto Frias, s/n                         \\
        4200-464 Porto, Portugal                                     \\
        jmartins@fe.up.pt  }
\date{}

\begin{document}
%%%%%%%%%%%%%%%%%%%%%1

\maketitle

\begin{abstract}
In this report a lumped transfer function model  for High Pressure Natural Gas Pipelines is derived. Starting with a partial nonlinear differential equation (PDE) model a high order  continuous state space (SS) linear model is obtained using a finite difference method. Next, from the SS representation an infinite order transfer function (TF) model is calculated.  In the end, this TF is approximated by a compact non-rational function.
\end{abstract}

\pagenumbering{arabic} \setcounter{page}{1}

\section{Introduction}
In this report we investigate the problem of the representation of a high pressure gas pipeline by a compact non rational transfer function model.
This model is used to simulate mass flow and pressure in a small high pressure pipeline, and
although this is a simple model with few parameters,  it seems to  have an accuracy comparable to the SIMONE$^{\textregistered }$  simulator.  Since this kind of models are suitable to control design and are well understood by control practitioners, it is our intention to apply them to gas leakage detection and gas network control.

\section{STATE-SPACE DISCRETE-IN-SPACE MODEL}
The gas dynamics within the pipes is represented by a set of partial differential equations (PDE). If we neglect the viscous and the turbulent effects of the flow and  assume small temperature changes within the gas and small heat exchanges with the surroundings of the pipeline, it can be described by the one-dimensional hyperbolic model
\begin{equation}
  \label{1}
  \left\{
  \begin{array}{l}
  \small \displaystyle \dfrac{\partial q (\ell,t)}{\partial t} = - \mathcal{A} \dfrac{\partial  p(\ell,t)}{\partial \ell}-
  \displaystyle \dfrac{f_c c^2}{2 \mathcal{D}\mathcal{A}}\dfrac{ q^2(\ell,t)}{p(\ell,t)}\\[2mm]
  \small \displaystyle \dfrac{\partial p(\ell,t)}{\partial t} = - \dfrac{c^2}{\mathcal{A}} \dfrac{\partial
  q (\ell,t)}{\partial \ell},
  \end{array} \right.
\end{equation}
where $\ell$ is space, $t$ is time, $p$ is edge pressure-drop, $q$ is mass flow, $\mathcal{A}$ is
the cross-sectional area, $\mathcal{D}$ is the  pipe diameter, $c$ is the isothermal speed of sound, and
$f_c$ is the friction factor.

In this research we linearised model \eqref{1}  around the operational levels $\left( p_{m}(\ell),q_{m}\right),$ where we assume a constant flow rate, and from the first equation of  \eqref{1}
$$p_{m}(\ell) =\sqrt{p_{m}^{2}({\ell}_{0})- \dfrac{f_c c^2}{2 \mathcal{D}\mathcal{A}^{2}} q_{m}^{2} (\ell-{\ell}_{0})}.$$
Hence we set $p(\ell,t)=p_m(\ell)+\Delta p(\ell,t)$ and $q(\ell,t)=q_m+\Delta q(\ell,t),$ where $\Delta p(\ell,t)$  and $\Delta q(\ell,t)$ are deviations from the pressure/flow operational levels, respectively.  Then $\dfrac{q^{2}(\ell,t)}{p(\ell,t)}=\dfrac{\left( q_m+\Delta q(\ell,t) \right)^{2}}{p_m(\ell)+\Delta p(\ell,t)}=\dfrac{q_{m}^{2}}{p_{m}(\ell)}+2\dfrac{q_{m}}{p_{m}(\ell)} \Delta q(\ell,t)- \dfrac{q_{m}^{2}}{p_{m}^{2}(\ell)} \Delta p(\ell,t)$

The third term may be neglected since the distribution networks operate at very high pressure, ca. 80 bar.  Then we substitute the remaining in the first equation
\[  \dfrac{\partial q(\ell,t)}{\partial t} = - \mathcal{A} \dfrac{\partial p(\ell,t)}{\partial \ell}-
  \displaystyle  \dfrac{f_c c^2}{2 \mathcal{D}\mathcal{A}} \dfrac{q_{m}}{p_{m}}\left( q_{m} +2 \Delta q(\ell,t) \right) . \]
Assumming small oscillations, $ \Delta q(\ell,t)  \approx 2 \Delta q(\ell,t),$ we may have $\left( q_{m} +2 \Delta q(\ell,t) \right) \approx q(\ell,t)$ and obtain the following linearized model:
\begin{equation}\label{3}
  \left\{
  \begin{array}{l}
     \dfrac{\partial q(\ell,t)}{\partial t} = - \mathcal{A} \dfrac{\partial p(\ell,t)}{\partial \ell}-2\alpha q(\ell,t)\\[2mm]
     \dfrac{\partial p(\ell,t)}{\partial t} = - \dfrac{c^2}{\mathcal{A}} \dfrac{\partial q(\ell,t)}{\partial \ell}.
\end{array} \right.
\end{equation}
where
\begin{equation}\label{2}
\alpha=  \dfrac{f_c c^2}{4 \mathcal{D}\mathcal{A}} \dfrac{q_{m}}{p_{m}}.
\end{equation}
Next, decompose the pipeline into sections $\mathcal{L}_i=[{\ell }_{i-1}, {\ell }_i]$, $i=1,2,\dots,N$,
where ${\ell }_0=0$, ${\ell }_N=L$ and $L$ is the length of the pipeline. We assume the massflow to be the same in each section and  accordingly define the following notation:
\begin{equation}\label{4}
  \begin{array}{rcl}
     q_0(t)&=&q(0,t)\\
     q_i(t)&=&q(\ell,t),\quad {\ell }_{i-1} < {\ell } < {\ell }_i,\quad i=1,2,\dots,N \\
     q_{N+1}(t)&=&q(L,t)\\
     p_i(t)&=&p(\ell_i,t), \quad i=0,1,\dots,N.
  \end{array}
\end{equation}
Making
\begin{equation}\label{5}
   \left.\dfrac{\partial \cdot(\ell,t)}{\partial \ell}\right|_{{\ell }={\ell }_i}\approx
\dfrac{\cdot(\ell_i,t)-\cdot(\ell_{i-1},t)}{\ell_i-\ell_{i-1}},\quad i=1,2,\dots,N,
\end{equation}
we can now approximate the linearized PDE \eqref{3} by the following discrete-in-space model
\begin{equation}\label{6}
   \begin{array}{rcl}
       \dot{q}_i(t)&=&\hspace{-2mm}\dfrac{\mathcal{A}}{\Delta \ell}\left[p_{i-1}(t)-p_i(t)\right]-2\alpha q_i(t),\\
       & &\hspace{25mm}i=1,\dots,N\\
       \dot{p}_{j-1}(t)&=&\hspace{-2mm}\dfrac{c^2}{\mathcal{A} \Delta \ell}\left[q_{j-1}(t)-q_j(t)\right]\\
       & &\hspace{25mm}j=1,\dots,N+1.
   \end{array}
\end{equation}
where
\begin{equation}\label{7}
   \Delta {\ell }={\ell }_{i+1}-{\ell }_{i}=\dfrac{L}{N},\quad i=0,\dots,N-1.
\end{equation}
The pipe can then be  described by the following state-space model:
\begin{equation}\label{8}
   \begin{array}{rcl}
      \dot{x}_1(t)&=&-\dfrac{c^2}{\mathcal{A} \Delta \ell}x_{N+2}(t)+\dfrac{c^2}{\mathcal{A} \Delta \ell}u_1(t)\\[2mm]
      \dot{x}_i(t)&=&\dfrac{c^2}{\mathcal{A} \Delta \ell}x_{N+i}(t)-\dfrac{c^2}{\mathcal{A} \Delta \ell}x_{N+i+1}(t)\\
      & &\\
      \dot{x}_{N+1}(t)&=&\dfrac{c^2}{\mathcal{A} \Delta \ell}x_{2N+1}(t)-\dfrac{c^2}{\mathcal{A} \Delta \ell}u_2(t)\\[2mm]
      \dot{x}_{N+1+j}(t)&=&\dfrac{\mathcal{A}}{\Delta \ell}x_{j}(t)-\dfrac{\mathcal{A}}{\Delta \ell}x_{j+1}(t) -2\alpha x_{N+1+j}(t)\\[4mm]
      y_1(t)&=&x_1(t)\\[2mm]
      y_2(t)&=&x_{N+1}(t),
   \end{array}
\end{equation}
where $i=1,\dots,N,$  $j=1,\dots N$ and also
\begin{equation}\label{9}
   \begin{array}{l}
      u(t)=\left[\ba{cc} q_0(t) &q_{N+1}(t)\ea\right]^T=\left[\ba{cc} u_1(t) & u_2(t)\ea\right]^T\\[2mm]
      x(t)=\left[\ba{ccccccc} p_0(t)  & \cdots & p_N(t) & | & q_1(t)  & \cdots & q_N(t) \ea\right]^T\\[2mm]
     y(t)=\left[\ba{cc} p_0(t) & p_{N}(t) \ea\right]^T=\left[\ba{cc} y_1(t) & y_2(t) \ea\right]^T.
   \end{array}
\end{equation}
In matrix notation:
\begin{equation}\label{10}
  \begin{array}{rcl}
     \dot{x}(t) & = & Ax(t)+Bu(t)\\
     y(t) & = & Cx(t).
     \end{array}
\end{equation}
Partition $A$ as:
\begin{equation}\label{11b}
A = \left[\ba{c|c} A_{11} & A_{12} \\ \hline A_{21} & A_{22} \ea\right]
\end{equation}
where  \begin{eqnarray}
 \label{12}
 A_{11} &=& 0_{(N+1) \times (N+1)} \\[4mm]
 \label{13}
 A_{12} &=&
 \scriptsize
 \left[
 \ba{ccccc}
    -\dfrac{c^2}{\mathcal{A} \Delta \ell} & 0 & \cdots & 0 & 0\\
    \dfrac{c^2}{\mathcal{A} \Delta \ell}  & -\dfrac{c^2}{\mathcal{A} \Delta \ell} & \cdots  & 0  & 0 \\
%    \vdots & \ddots &  \ddots & \cdots  & \cdots\\
    \vdots & \vdots &  \ddots &  \vdots & \vdots\\
    0  & 0 & \cdots  & \dfrac{c^2}{\mathcal{A} \Delta \ell} & -\dfrac{c^2}{\mathcal{A} \Delta \ell}\\
    0  & 0 & \cdots  & 0 & \dfrac{c^2}{\mathcal{A} \Delta \ell}
 \ea
 \right] \in \mathbb{R}{(N+1)\times N} \\[4mm]
 \label{14}
 A_{21}&=&
 \left[
 \ba{ccccc}
    \dfrac{\mathcal{A}}{\Delta \ell} & - \dfrac{\mathcal{A}}{\Delta \ell} & \cdots & 0 & 0 \\
     \vdots &  \vdots & \ddots & \vdots & \vdots\\
    0 & 0 &  \cdots & \dfrac{\mathcal{A}}{\Delta \ell} & - \dfrac{\mathcal{A}}{\Delta \ell}
 \ea
 \right]  \in \mathbb{R}^{N\times (N+1) }\\[4mm]
 \label{15}
 A_{22}&=&-2\alpha I_N\\[4mm]
 \label{16}
 B&=& \dfrac{c^2}{\mathcal{A} \Delta \ell} \left[\ba{c|c} e_1 & -e_{N+1} \ea \right]\\[4mm]
 \label{17}
 C&=&\left[\ba{c|c} e_1 & e_{N+1} \ea \right]^T
\end{eqnarray}
where $e_i$ is the $i^{th}$ vector of the canonical orthonormal basis, i.e., a vector with the
the $i^{th}$ component equal to one and the others equal to zero.
\section{Spectral analysis of A}
In order to learn more about the system \eqref{2}, we analyse the spectrum of matrix A.
Therefore the following theorem:
\begin{theorem}\label{Th1}
The eigenvalues of $A$ defined in \eqref{12}--\eqref{15} are
\begin{eqnarray}
\lambda_0&=&0\\
\lambda_{\pm k}&=&
-\dfrac{f_c c^2 Q_{m}}{4 \mathcal{D} \mathcal{A} P_{m}} \pm j
\sqrt{\left(2\dfrac{c}{\Delta \ell}\sin\left(\frac{k\pi}{2(N+1)}\right)\right)^2
-\left(\dfrac{f_c c^2 Q_{m}}{4 \mathcal{D} \mathcal{A} P_{m}}\right)^2},\;k= 1,\dots,  N \label{32}
\end{eqnarray}
\end{theorem}
\textbf{Proof:}
In the Appendix~\ref{changeofvariables} it was proven that matrix $\bar{A}$
is equivalent to $A$ up to a similarity transformation. Consequently they have the same eigenvalues and using \eqref{A-4}, we have:
\bea
\det(sI_{2N+1}-A)=\det(sI_{2N+1}-\bar{A})=
s\det\left[\ba{c|c} sI_N-\bar{A}_{11} & -\bar{A}_{12}  \\ \hline \\[-4mm]-\bar{A}_{21} & sI_N-\bar{A}_{22} \ea \right]=0
\eea
And this is equivalent to $s=0$ and $\det\left[\ba{c|c} sI_N-\bar{A}_{11} & -\bar{A}_{12}  \\ \hline\\[-4mm] -\bar{A}_{21} & sI_N-\bar{A}_{22} \ea \right]=0.$
Therefore,  $A$ has a zero eigenvalue, that is, $\lambda_0=0.$

From Fact 2.13.10 in \cite[pp. 62--63]{Bernstein}, we have that  for arbitrary matrices
$\mathbf{A},\;\mathbf{B},\;\mathbf{C}$ and $\mathbf{D} \in \cR^{N \times N}$ such that
$\mathbf{A}\mathbf{B}=\mathbf{B}\mathbf{A}$ then
\bea
\det\left[\ba{c|c} \mathbf{A} & \mathbf{B} \\\hline \mathbf{C} & \mathbf{D} \ea \right]=
\det(\mathbf{D}\mathbf{A}-\mathbf{C}\mathbf{B}).
\eea
Thus, if we take
\bea
\mathbf{A}&=&sI_N-\bar{A}_{11}\\
\mathbf{B}&=&-\bar{A}_{12}\\
\mathbf{C}&=&-\bar{A}_{21}\\
\mathbf{D}&=&sI_N-\bar{A}_{22}
\eea
we see that $(sI_N-\bar{A}_{11})(-\bar{A}_{12})=(-\bar{A}_{12})(sI_N-\bar{A}_{11})\Rightarrow \mathbf{A}\mathbf{B}=\mathbf{B}\mathbf{A}$ because
$sI_N-\bar{A}_{11}$ is a diagonal matrix. Consequently
\bea
\det\left[\ba{c|c} sI_N-\bar{A}_{11} & -\bar{A}_{12} \\\hline\\ [-4mm] -\bar{A}_{21} & sI_N-\bar{A}_{22} \ea \right]
=\det\bigl((sI_N-\bar{A}_{22})(sI_N-\bar{A}_{11})-\bar{A}_{21}\bar{A}_{12}\bigr).
\eea
Given that $\bar{A}_{11}=0_{N\times N}$ and $\bar{A}_{22}=-2\alpha I_N,$ then
\bea
&&\det\left[\ba{c|c} sI_N-\bar{A}_{11} & -\bar{A}_{12} \\\hline\\[-4mm] -\bar{A}_{21} & sI_N-\bar{A}_{22} \ea \right]=
\det\left((s^2+2\alpha s)I_N-\bar{A}_{21}\bar{A}_{12}\right)\\
&=&
\det\left((s^2+2\alpha s+\alpha^2)I_N-\bar{A}_{21}\bar{A}_{12}-\alpha^2I_N\right)
=\det\left((s+\alpha)^2I_N-\bar{A}_{21}\bar{A}_{12}-\alpha^2I_N\right).
\eea
If we define the following the change of variable:
\bean
\label{34A}
\mathcal{S}=(s+\alpha)^2
\eean
then we can write
\bean
\label{35A}
\det\left[\ba{c|c} sI_N-\bar{A}_{11} & -\bar{A}_{12} \\\hline\\[-4mm] -\bar{A}_{21} & sI_N-\bar{A}_{22} \ea \right]&=&\det\left(\mathcal{S}I_N-\left(\bar{A}_{21}\bar{A}_{12}+\alpha^2I_N\right)\right)
\eean
From this equation, the eigenvalues of $\bar{A}_{21}\bar{A}_{12}+\alpha^2I_N$ that we denote by $\Lambda\left(\bar{A}_{21}\bar{A}_{12}+\alpha^2 I_N\right)$,
are the values of $\mathcal{S}=(s+\alpha)^2$ that also set $\det(sI_N-A)$ to zero.  From \eqref{35A}, $\Lambda\left( \bar{A}_{21}\bar{A}_{12}+ \alpha^2 I_N \right)=\left(\Lambda(A)+\alpha\right)^2$, where $\Lambda(A)$ denotes the non zero eigenvalues of $A$.

From the eigenvalues properties,
\bean
\label{36A}
\Lambda\left(\bar{A}_{21}\bar{A}_{12}+\alpha^2I_N\right)=\Lambda\left(\bar{A}_{21}\bar{A}_{12}\right)+\alpha^2.
\eean
The product $\bar{A}_{21}\bar{A}_{12}$ is
\bea
\bar{A}_{21}\bar{A}_{12}=\left(\dfrac{c}{\Delta \ell}\right)^2
\left[\ba{cccccc} -2 & 1 & 0 & \cdots & 0 & 0 \\ 1 & -2 & 1 & \cdots & 0 & 0 \\ 0 & 1 & -2 & \ddots & 0 & 0 \\
\vdots & \vdots & \ddots & \ddots & \ddots & \vdots \\ 0 & 0 & 0 & \ddots & -2 & 1 \\ 0 & 0 & 0 & \cdots & 1 & -2\ea\right].
\eea

Using Fact 5.10.25 in \cite[pp. 200]{Bernstein}
\bean
\label{37A}
\Lambda\left(\bar{A}_{21}\bar{A}_{12}\right)=
-2\left(\dfrac{c}{\Delta \ell}\right)^2\left( 1-\cos\left(\frac{k\pi}{N+1}\right)\right),\;k=1,\dots,N
\eean
Then
\bean
\label{38A}
\Lambda\left(\bar{A}_{21}\bar{A}_{12}+\alpha^2I_N\right)=
-2 \left(\dfrac{c}{\Delta \ell}\right)^2\left(1-\cos\left(\frac{k\pi}{N+1}\right)\right)+\alpha^2,\;k=1,\dots,N
\eean
are the values of $\mathcal{S}=\left(s+\alpha\right)^2$ that set the characteristic equation of $A$ to zero.

Consequently
\bea
\left(\Lambda(A)+\alpha\right)^2=\Lambda\left(\bar{A}_{21}\bar{A}_{12}+\alpha^2I_N\right)
\eea
and from \eqref{34A}--\eqref{38A} the eigenvalues of $A$ are
$ \Lambda(A)=-\alpha\pm\sqrt{
-2\left(\dfrac{c}{\Delta \ell}\right)^2\left(1-\cos\left(\dfrac{k\pi}{N+1}\right)\right)+\alpha^2}$\\
That is
\bea
\Lambda(A)&=&-\alpha \pm j \sqrt{2\left(\dfrac{c}{\Delta \ell}\right)^2\left(1-\cos\left(\frac{k\pi}{N+1}\right)\right)-\alpha^2}\\[4mm]
&=&-\alpha  \pm j
\sqrt{4\left(\dfrac{c}{\Delta \ell}\right)^2\sin^2\left(\dfrac{k\pi}{2(N+1)}\right)
-\alpha^2}
\eea
Recalling the definition of $\alpha $ in equation \eqref{2}, we have:
\begin{eqnarray*}
\Lambda(A)&=&-\dfrac{f_c c^2 Q_{m}}{4 \mathcal{D} \mathcal{A} P_{m}} \mp j
\sqrt{\left(2\dfrac{c}{\Delta \ell}\sin\left(\frac{k\pi}{2(N+1)}\right)\right)^2
-\left(\dfrac{f_c c^2 Q_{m}}{4 \mathcal{D} \mathcal{A} P_{m}}\right)^2},\;k= 1,\dots, N \label{39}
\end{eqnarray*}
and this completes the proof.

\hfill $\Box$

The asymptotic case of the nonzero eigenvalues is reported in the following corollary
\begin{corollary}\label{Cor1}
If $N\rightarrow\infty$ then the eigenvalues of $A$ are
\begin{eqnarray}
\lambda_0&=&0\\
\lambda_{\pm k}
&=&-\underbrace{\dfrac{f_c c^2 Q_{m}}{4 \mathcal{D} \mathcal{A} P_{m}} }_{:= \alpha}\mp j
\underbrace{\sqrt{\left(\dfrac{k\pi}{T_d}\right)^2-\left(\dfrac{f_c c^2 Q_{m}}{4 \mathcal{D} \mathcal{A} P_{m}}\right)^2}}_{:=b_{k}}
,\;k= 1, 2,\dots\label{assymeigen}
\end{eqnarray}
where $L$ is the pipe length and $T_d$ the time that a mass pressure takes to cross the pipeline,  between its boundaries, at a constant speed $c$.
\end{corollary}
\textbf{Proof:} $\displaystyle \lim_{N\rightarrow \infty}\lambda_0=0$ is trivial\\
Since $\Delta \ell$ is given by
\bea
\Delta \ell=\dfrac{L}{N}
\eea
then
\bea
\dfrac{c}{\Delta \ell}\sin\left(\frac{k\pi}{2(N+1)}\right)=\dfrac{cN}{L}\sin\left(\frac{k\pi}{2(N+1)}\right).
\eea
Taking the limit when $N\rightarrow\infty$, $\displaystyle\lim_{N\rightarrow \infty}\dfrac{cN}{L}\sin\left(\frac{k\pi}{2(N+1)}\right)=\dfrac{c}{2L}k\pi $

and, consequently,
\bea
\lim_{N\rightarrow\infty}\lambda_{\pm k}&=&-\dfrac{f_c c^2 Q_{m}}{4 \mathcal{D} \mathcal{A} P_{m}} \mp j
\sqrt{\left(2\dfrac{c}{2L}k\pi\right)^2-\left(\dfrac{f_c c^2 Q_{m}}{4 \mathcal{D} \mathcal{A} P_{m}}\right)^2}\\
&=&-\dfrac{f_c c^2 Q_{m}}{4 \mathcal{D} \mathcal{A} P_{m}} \mp j
\sqrt{\left(\dfrac{c}{L}k\pi\right)^2-\left(\dfrac{f_c c^2 Q_{m}}{4 \mathcal{D} \mathcal{A} P_{m}}\right)^2},\;k= 1, 2,\dots
\eea
Given that $T_d=\dfrac{L}{c}$, one obtains
\bea
\lim_{N\rightarrow\infty}\lambda_{\mp k}=-\dfrac{f_c c^2 Q_{m}}{4 \mathcal{D} \mathcal{A} P_{m}} \mp j
\sqrt{\left(\dfrac{k\pi}{T_d}\right)^2-\left(\dfrac{f_c c^2 Q_{m}}{4 \mathcal{D} \mathcal{A} P_{m}}\right)^2}
,\;k= 1, 2,\dots
\eea
and this completes the proof.

\hfill$\Box$

The zero eigenvalue means that there is an integrator in the pipeline model.

It has associated an eigenvector $v_0$, which defines a direction in the state-space where the pipeline behaves like a pure integrator.   The following lemma gives the value of this eigenvector:

\begin{lemma}\label{L1}
Consider $A$ as in \eqref{12}--\eqref{15}.
Then
\bean
\label{26}
v_0&=&e_1+e_2+\dots+e_{N+1}
\eean
where $e_i$ is the $i^{th}$ vector of the canonical orthonormal base in $\cR^{2N+1}$, is the eigenvector associated  to the zero eigenvalue, i. e.,
$Av_0=0.$
\end{lemma}

\textbf{Proof:}
Let us denote $v_0$ as
\bean
\label{28}
v_0&=& \left[\ba{c} v \\ 0_{N}\ea\right]
\eean
where $0_N$ is the zero vector in $R^N$ and
\bean
\label{29A}
v=\varepsilon_1+\varepsilon_2+\dots+\varepsilon_{N+1} \in \cR^{N+1}
\eean
where $\varepsilon_i$ is the $i^{th}$ vector of the canonical orthonormal basis in $R^{N+1}$.
Using \eqref{11b} and \eqref{28}
\bean
\label{30}
Av_0&=& \left[\ba{c} A_{11}v \\\hline A_{21} v \ea\right]
\eean
Given that $A_{11}=0_{(N+1) \times (N+1)}$ then $A_{11}v=0$. From \eqref{14} and \eqref{29A}
\bea
A_{21}v&=&\left[\ba{cc} \dfrac{\mathcal{A}}{\Delta \ell}-\dfrac{\mathcal{A}}{\Delta \ell}\\[2mm]
                        \dfrac{\mathcal{A}}{\Delta \ell}-\dfrac{\mathcal{A}}{\Delta \ell}\\[2mm]
                        \vdots\\
                        \dfrac{\mathcal{A}}{\Delta \ell}-\dfrac{\mathcal{A}}{\Delta \ell}
           \ea \right]=0_{N}
\eea
and this completes the proof.

\hfill$\Box$

\begin{remark}
In Lemma~\ref{L1}, to prove the existence of the eigenvalue associated to the zero eigenvalue we only used the submatrices $A_{11}$ and $A_{12}$. Therefore, we can say that this eigenvalue is generated by these submatrices. The nonlinearity of the model is only expressed by matrix $A_{22}.$  For this reason, if we decompose the full nonlinear model into a linear subsystem
in cascade with a nonlinear one, the zero eigenvalue would appear in the linear subsystem indicating the presence
of an integrator in the full  model.    Also, $A_{11}$ and $A_{12}$
depend neither on $p_{m}$ nor  on $q_{m}$.
\end{remark}
%==================================
\section{Transfer functions characterisation}
We determine the transfer function.  Recall that the massflow at the boundaries were chosen to be our inputs and the pressure at the boundaries  our outputs.

To start, we apply the Laplace transform to \eqref{10} and  obtain:
\[
Y(s)=C\left( sI-A \right)^{-1}B U(s)
\]
where $Y{(s)}=\left[\ba{cc} Y_1(s)&
Y_2(s) \ea \right]^T$ and $U(s)=\left[\ba{cc} U_1(s) & U_2(s) \ea\right]^T,$ and $F(s)$ denotes the Laplace transform of $f(t)$.  That is
\[
\left[ \begin{array}{c}
Y_1(s) \\ Y_2(s)
\end{array}
\right]= \left[ \begin{array}{cc}
G_{11}(s) & G_{12}(s)\\
G_{21}(s) & G_{22}(s)
\end{array}
\right] \left[ \begin{array}{c}
U_1(s) \\ U_2(s)
\end{array}
\right] .
\]
Also:
\begin{eqnarray}
 \left[ \begin{array}{c}
Y_1(s)\\ Y_2 (s)
\end{array}
\right]&=& \left[ \begin{array}{c}
C_1\\
C_2
\end{array}
\right] \left( sI-A \right)^{-1}\left[ \begin{array}{cc}
B_1 & B_2
\end{array}
\right]  \left[ \begin{array}{c}
U_1(s) \\ U_2(s)
\end{array}
\right] \nonumber \\[2mm]
 \left[ \begin{array}{c}
Y_1(s) \\ Y_2(s)
\end{array}
\right]&=& \left[ \begin{array}{cc}
C_1 \left( sI-A \right)^{-1} B_1&C_1 \left( sI-A \right)^{-1}B_2 \\
C_2 \left( sI-A \right)^{-1} B_1& C_2 \left( sI-A \right)^{-1}B_2
\end{array}
\right] \left[ \begin{array}{c}
U_1(s) \\ U_2(s)
\end{array}
\right] \label{TFcomponents}
\end{eqnarray}
and also \bean
\label{40}
\ba{rcl}
G_{11}(s)&=&={\left.  \dfrac{Y_1(s)}{U_1(s)}\right|}_{U_{2}(s)=0}={\left. \dfrac{P_{0}(s)}{Q_{0}(s)}\right|}_{Q_{N+1}(s)=0}\\[8mm]
G_{22}(s)&=&{\left. \dfrac{Y_2(s)}{U_2(s)}\right|}_{U_{1}(s)=0}={\left. \dfrac{P_{N}(s)}{Q_{N+1}(s)}\right|}_{Q_{0}(s)=0}\\[8mm]
G_{12}(s)&=&{\left. \dfrac{Y_1(s)}{U_2(s)}\right|}_{U_{1}(s)=0}={\left. \dfrac{P_{0}(s)}{Q_{N+1}(s)}\right|}_{Q_{0}(s)=0}\\[8mm]
G_{21}(s)&=&{\left. \dfrac{Y_2(s)}{U_1(s)}\right|}_{U_{2}(s)=0}={\left. \dfrac{P_{N}(s)}{Q_{0}(s)}\right|}_{Q_{N+1}(s)=0}.
\ea
\eean
\subsection{Transfer function $G_{11}$}
When we select this transfer function, it means that we are interested in the transfer function between
the pressure and massflow at the intake node, i.e. $\dfrac{Y_1(s)}{U_1(s)}$ when $U_{2}(s)=0$   that is:
$$G_{11}(s)={\left.\dfrac{P_0(s)}{Q_0(s)}\right|_{Q_1(s)=0}}={\left.\dfrac{Y_1(s)}{U_1(s)}\right|_{U_2(s)=0}}$$
and hence:
\bean
\label{43}
G_{11}(s)&=&C_1\left(sI-A\right)^{-1}B_1,
\eean
where $B_1$ and $C_1$ are the first column and first row of $B$ and $C$ in \eqref{16}--\eqref{17}, respectively.  $G_{11}(s)$ is a rational function
whose poles are the eigenvalues of $A$.

The following theorem states the zeros of this transfer function.
\begin{theorem}\label{Th2}
The zeros of $G_{11}(s)$ are
\begin{eqnarray}
z_{ \pm k}&=&
-\dfrac{f_c c^2 Q_{m}}{4 \mathcal{D} \mathcal{A} P_{m}} \pm j
\sqrt{\left(2\dfrac{c}{\Delta \ell}\sin\left(\dfrac{(2k-1)\pi}{2(2N+1)}\right)\right)^2
-\left(\dfrac{f_c c^2 Q_{m}}{4 \mathcal{D} \mathcal{A} P_{m}}\right)^2},\quad k= 1,\dots, N\label{48}
\end{eqnarray}
\end{theorem}
\textbf{Proof:}
To do this we recall the result from \cite[pp. 284]{Carvalho:1993}, we have that the zeros of the transfer function \eqref{43} are the zeros of the following polynomial
\begin{equation}
\left|\begin{array}{c|c}
sI_{(2N+1)} -A & -B_1 \\\hline
C_1 & 0
\end{array}
\right|=0 .
\end{equation}
Recalling that $C_1=e_1 \in \cR^{(2N+1)}$ and $B_1=\dfrac{c^2}{{\cal A} \Delta \ell }e_1\cR^{(2N+1)}$ (see equations \eqref{16}--\eqref{17}), we have:
\begin{eqnarray}
\left|\begin{array}{c|c}
sI_{(2N+1)} -A & -e_1\\\hline
e_1^T& 0
\end{array}
\right|&=& 0 \quad \Leftrightarrow\\[2mm]
\Leftrightarrow \left|\begin{array}{c|c}
\left[ \begin{array}{c|c} sI_{N+1}-A_{11} & -A_{12} \\\hline  -A_{21} & sI_N-A_{22}\end{array}\right] & \begin{array}{c} -1 \\ 0 \\  \vdots \\ 0 \end{array}\\[2mm]\hline
\begin{array}{ccccccc}1 & 0& \cdots & \cdots& \cdots& \cdots & 0\end{array}&  0
\end{array}
\right|
&=& 0\quad \Leftrightarrow
\end{eqnarray}
using the definition of $A_{11}$ and $A_{22}$ in \eqref{12} and \eqref{15}
\begin{eqnarray}
\Leftrightarrow\left|\begin{array}{c|c}
\left[ \begin{array}{c|c} sI_{N+1} & -A_{12} \\\hline -A_{21} & (s+2\alpha)I_N\end{array}\right] & \begin{array}{c} -1 \\ 0 \\  \vdots \\ 0 \end{array}\\[2mm]\hline
\begin{array}{cccccc}1 & 0 & \cdots& \cdots& \cdots & 0\end{array}&  0
\end{array}
\right|&=& 0 .
\end{eqnarray}
We develop this determinant first along the last column and next along the last row, and obtain:
\begin{eqnarray}
\left| \begin{array}{c|c} sI_{N} & -\bar{A}_{12} \\[4mm] \hline \\ -\bar{A}_{21} & (s+2\alpha)I_N\end{array}
\right|&=& 0
\end{eqnarray}
where $\bar{A}_{12}$ denotes matrix $A_{12}$ without the $1$st row and  $\bar{A}_{21}$ denotes matrix $A_{21}$ without the $1$st column.
Next, as $sI_{N}\left(-\bar{A}_{12} \right)=\left(-\bar{A}_{12} \right)sI_{N} $, we apply again
Fact 2.13.10 in \cite[pp. 62--63]{Bernstein}, which states that for the arbitrary matrices
$\mathbf{A},\;\mathbf{B},\;\mathbf{C}$ and $\mathbf{D} \in \cR^{N \times N}$ such that
$\mathbf{A}\mathbf{B}=\mathbf{B}\mathbf{A}$ then
\bea
\det\left[\ba{c|c} \mathbf{A} & \mathbf{B} \\\hline \mathbf{C} & \mathbf{D} \ea \right]=
\det(\mathbf{D}\mathbf{A}-\mathbf{C}\mathbf{B})
\eea
and
\begin{eqnarray}
\left| \begin{array}{c|c} sI_{N} & -\bar{A}_{12} \\ \hline  -\bar{A}_{21} & (s+2\alpha)I_N\end{array}
\right|&=& \left| s(s+2\alpha)I_N-\bar{A}_{21}\bar{A}_{12}\right|\\[2mm]
&=& \left|\left(s^2+2\alpha s+\alpha^2 \right) I_N-\left( \bar{A}_{21}\bar{A}_{12} +\alpha^2\right) \right|\\[2mm]
&=&\left|\left(s+\alpha\right)^2  I_N-\left( \bar{A}_{21}\bar{A}_{12} +\alpha^2 \right) \right|
\end{eqnarray}
We do the same change of variable as before
\begin{equation}\label{chvar}
\mathcal{S}=(s+\alpha)^2
\end{equation}
and then can write:
\begin{eqnarray}
\left| \begin{array}{c|c} sI_{N} & -\bar{A}_{12} \\ \hline  -\bar{A}_{21} & (s+2\alpha)I_N\end{array}
\right|
&=&\left|  \mathcal{S}I_N-\left( \bar{A}_{21}\bar{A}_{12} +\alpha^2 \right) \right|
\end{eqnarray}
Next, we calculate the spectrum of matrix $\left( \bar{A}_{21}\bar{A}_{12} + \alpha^2 \right),$ that is:
\[
\Lambda \left( \bar{A}_{21}\bar{A}_{12} +\alpha^2\right)=\Lambda \left( \bar{A}_{21}\bar{A}_{12} \right)+\alpha^2
\]
Now, we calculate the product:
\begin{eqnarray*}
\bar{A}_{21}\bar{A}_{12}& =&  \left[
\ba{ccccc}
- \dfrac{\mathcal{A}}{\Delta \ell} & 0 &\cdots & 0 & 0 \\
[4mm]
\dfrac{\mathcal{A}}{\Delta \ell} & - \dfrac{\mathcal{A}}{\Delta \ell} & \cdots & 0 & 0\\
0 & \ddots & \ddots & \cdots & \cdots\\
\vdots  & \cdots & \ddots & \ddots & \cdots\\
0 & 0 & \cdots & \dfrac{\mathcal{A}}{\Delta \ell} & - \dfrac{\mathcal{A}}{\Delta \ell}
\ea
\right]\left[
\ba{ccccccc}
\dfrac{c^2}{\mathcal{A} \Delta \ell} & -\dfrac{c^2}{\mathcal{A} \Delta \ell} & 0 & 0 & \cdots & 0 & 0\\
[4mm]
0 &\dfrac{c^2}{\mathcal{A} \Delta \ell} & -\dfrac{c^2}{\mathcal{A} \Delta \ell} & 0 & \cdots & 0 & 0 \\
[4mm]
\vdots & \vdots & \ddots & \ddots &  \cdots & \vdots & \vdots\\
0 &0&0&0 &\cdots &\dfrac{c^2}{\mathcal{A} \Delta \ell} & -\dfrac{c^2}{\mathcal{A} \Delta \ell}  \\
0 &0 & 0 & 0 & \cdots & 0 &\dfrac{c^2}{\mathcal{A} \Delta \ell}
\ea
\right]\\[4mm]
& =& (\dfrac{c}{\Delta \ell})^2  \left[\ba{cccccc} -1 & 1 & 0 & \cdots & 0 & 0 \\ 1 & -2 & 1 & \cdots & 0 & 0 \\ 0 & 1 & -2 & \ddots & 0 & 0 \\
\vdots & \vdots & \ddots & \ddots & \ddots & \vdots \\ 0 & 0 & 0 & \ddots & -2 & 1 \\ 0 & 0 & 0 & \cdots & 1 & -2\ea\right]\in \cR^N
\end{eqnarray*}
Then
\[ \Lambda\left(\bar{A}_{21}\bar{A}_{12}\right)= (\dfrac{c}{\Delta \ell})^2
\Lambda\left(M_N\right)
\]
where
\begin{equation}
\label{61}
M_N=\left[\ba{cccccc} -1 & 1 & 0 & \cdots & 0 & 0 \\ 1 & -2 & 1 & \cdots & 0 & 0 \\ 0 & 1 & -2 & \ddots & 0 & 0 \\
\vdots & \vdots & \ddots & \ddots & \ddots & \vdots \\ 0 & 0 & 0 & \ddots & -2 & 1 \\ 0 & 0 & 0 & \cdots & 1 & -2\ea\right]
\in \cR^{N \times N}.
\end{equation}
From \cite[pp. 72]{Yueh:2005}
\begin{equation}
\label{62}
\Lambda\left(M_N\right)=-2+2\cos\left(\dfrac{\left( 2k-1\right) \pi}{2N+1}\right),\quad k=1,2,3,\ldots,N.
\end{equation}
Having that
\begin{eqnarray}
\Lambda \left( \bar{A}_{21}\bar{A}_{12} +\alpha^2\right)& = & -2\left(\dfrac{c}{\Delta \ell}\right)^2\left( 1-\cos\left(\dfrac{\left( 2k-1\right) \pi}{2N+1}\right)  \right)+\alpha^2\\
& = & j^2 \left( 2\left(\dfrac{c}{\Delta \ell}\right)^2\left( 1-\cos\left(\dfrac{\left( 2k-1\right) \pi}{2N+1}\right)  \right)-\alpha^2\right) \label{63}
\end{eqnarray}
then from \eqref{chvar}:
\begin{equation*}
\label{64}
\ba{rcl}
z_{\pm k}&=&-\alpha \pm j \sqrt{2\left(\dfrac{c}{\Delta \ell}\right)^2\left(1-\cos\left(\dfrac{(2k-1)\pi}{2N+1}\right)\right)-\alpha^2}\\[4mm]
&=&-\dfrac{f_c c^2 Q_{m}}{4 \mathcal{D} \mathcal{A} P_{m}} \pm j
\sqrt{4\left(\dfrac{c}{\Delta \ell}\right)^2\sin^2\left(\dfrac{(2k-1)\pi}{2(2N+1)}\right)
-\left(\dfrac{f_c c^2 Q_{m}}{4 \mathcal{D} \mathcal{A} P_{m}}\right)^2}\\[4mm]
&=&-\dfrac{f_c c^2 Q_{m}}{4 \mathcal{D} \mathcal{A} P_{m}} \pm j
\sqrt{\left(2\dfrac{c}{\Delta \ell}\sin\left(\dfrac{(2k-1)\pi}{2\left( 2N+1\right) }\right)\right)^2
-\left(\dfrac{f_c c^2 Q_{m}}{4 \mathcal{D} \mathcal{A} P_{m}}\right)^2},\quad k= 1,\dots, N
\ea
\end{equation*}
and this completes the proof.

\hfill$\Box$

The following corollary resolves  the asymptotic case.
\begin{corollary}\label{Cor2}
If $N\rightarrow\infty$ then the zero of $G_{11}(s)$ are
\begin{equation}\label{zerosG11}
z_{ \pm k}=-\underbrace{\dfrac{f_c c^2 Q_{m}}{4 \mathcal{D} \mathcal{A} P_{m}}}_{:= \alpha } \pm j
\underbrace{\sqrt{\left(\dfrac{(2k-1)\pi}{2T_d}\right)^2-\left(\dfrac{f_c c^2 Q_{m}}{4 \mathcal{D} \mathcal{A} P_{m}}\right)^2}}_{:= \beta_{k}}
,\;k= 1, 2,\dots
\end{equation}
where $L$ is the pipe length and $T_d$ the time that a particle of gas  takes to cross the pipeline between its boundaries, at a constant speed
$c$.
\end{corollary}
\textbf{Proof:}
Since $\Delta \ell$ is given by
\bea
\Delta \ell=\dfrac{L}{N}
\eea
then
\bea
\dfrac{c}{\Delta \ell}\sin \left(\dfrac{(2k-1)\pi}{2\left( 2N+1\right) }\right)=\dfrac{cN}{L}\sin \left(\dfrac{(2k-1)\pi}{2\left( 2N+1\right) }\right).
\eea
Taking the limit when $N\rightarrow\infty$,
\bea
\lim_{N\rightarrow \infty}\dfrac{cN}{L}\sin \left(\dfrac{(2k-1)\pi}{2\left( 2N+1\right) }\right)=
\lim_{N\rightarrow \infty}\dfrac{cN \left(2k-1 \right) \pi}{2(2N+1)L}
\lim_{N\rightarrow \infty}\dfrac{\sin \left(
\dfrac{(2k-1)\pi}{2\left( 2N+1\right)}\right)}{\dfrac{\left(2k-1 \right) \pi}{2(2N+1)}}=\dfrac{c\left(2k-1 \right) \pi}{4 L}
\eea
and, consequently,
\bea
\lim_{N\rightarrow\infty}z_{\pm k}&=&-\dfrac{f_c c^2 Q_{m}}{4 \mathcal{D} \mathcal{A} P_{m}} \pm j
\sqrt{\left(\dfrac{c\left(2k-1 \right) \pi}{2 L}\right)^2-\left(\dfrac{f_c c^2 Q_{m}}{4 \mathcal{D} \mathcal{A} P_{m}}\right)^2}\quad k= 1, 2,\dots
\eea
Given that $T_d=\dfrac{L}{c}$ then
\bea
\lim_{N\rightarrow\infty}z_{\pm k}&=&-\dfrac{f_c c^2 Q_{m}}{4 \mathcal{D} \mathcal{A} P_{m}} \pm j
\sqrt{\left(\dfrac{\left(2k-1 \right) \pi}{2 Td}\right)^2-\left(\dfrac{f_c c^2 Q_{m}}{4 \mathcal{D} \mathcal{A} P_{m}}\right)^2}\quad k= 1, 2,\dots
\eea
and this completes the proof.

\hfill$\Box$

\begin{corollary}
\label{Cor3}
Transfer function $G_{11}$ has the following form:
\begin{equation}\label{TFG11}
G_{11}=\dfrac{K_{G} \displaystyle\prod_{k=1}^{\infty} \left( \dfrac{s^{2}}{z_{k}z_{-k}}+s\left( \dfrac{1}{z_{k}}+\dfrac{1}{z_{-k}}\right)+1 \right)}{s \displaystyle\prod_{k=1}^{\infty} \left( \dfrac{s^{2}}{\lambda_{k}\lambda_{-k}}+s \left(\dfrac{1}{\lambda_{k}}+\dfrac{1}{\lambda_{-k}}\right)+1 \right) }
\end{equation}
where $z_{k}=\alpha+j \beta_{k}$ and  $z_{-k}=\alpha-j \beta_{k},$ as well as $\lambda_{k}=\alpha+j b_{k}$ and  $\lambda_{-k}=\alpha-j b_{k},$ as defined in Corollary~\ref{Cor2} and  $\lambda_{i}$ are defined in Corollary~\ref{Cor1}.
\end{corollary}
\textbf{Proof:}
From Corollary~\ref{Cor2} and Corollary~\ref{Cor1}, we can write:
\begin{equation}\label{TFG11b}
G_{11}=\dfrac{K_{G} \displaystyle\prod_{k=1}^{\infty} \left( \dfrac{s}{-z_{k}}+1 \right) e^{\frac{s}{\beta_{k}}} \left(  \dfrac{s}{-z_{-k}}+1 \right)e^{\frac{s}{-\beta_{k}}}}{s \displaystyle\prod_{k=1}^{\infty} \left( \dfrac{s}{-\lambda_{k}}+1\right)e^{\frac{s}{b_{k}}} \left( \dfrac{s}{-\lambda_{-k}}+1 \right)e^{\frac{s}{-b_{k}}} }
\end{equation}
and expression \eqref{TFG11} follows immediately, after  calculating the products:

$ \left( \dfrac{s}{-z_{k}}+1 \right)  \left(  \dfrac{s}{-z_{-k}}+1 \right)  $ and
$ \left( \dfrac{s}{-\lambda_{k}}+1\right) \left( \dfrac{s}{-\lambda_{-k}}+1 \right) .$

\hfill$\Box$

To complete the transfer function characterisation we need to compute the gain $K_G$.

\begin{theorem}\label{Th3}
Consider  $B_{1},$ the first column of $B$ defined in \eqref{16}, written in the base $$\left\{ v_{-N},\ldots, v_{-1},v_{0},v_{1},\ldots, v_{N}\right\}$$ where  $v_i$, $i=-N,\dots,-1,0,1,\ldots,N,$ are the eigenvectores of A.  If $\vartheta_{0}$ is the component of $B_{1}$  along $v_{0}$ then
$$K_G=\vartheta_{0}.$$
\end{theorem}
\textbf{Proof:}
Denote
the zero eigenvalue of $A$ as $\lambda_0$ and $\lambda_i$, $i=-N,\dots-1,1,\dots,N$ the complex eigenvalues, i.e.
$\lambda_{-i}=\lambda_i^*$ where $\cdot^*$ means the conjugate eigenvalue. Since all eigenvalues have
multiplicity one there are 2N+1 independent eigenvectors $v_i$, $i=-N,\dots,N,$ respectively associated to each eigenvalue
$\lambda_i$. Thus $v_i$ and $(s-\lambda_i)$, $i=-N,\dots, -1,0,1,\dots,N$ are, respectively,
the eigenvectors and the eigenvalues of $sI-A$. Consequently,
\bean
\label{66}
\left(SI-A\right)^{-1}v_i&=&\dfrac{1}{s-\lambda_i}v_i
\eean
i.e., $\dfrac{1}{s-\lambda_i}$, $i=-N,\;\dots,\;N$ are eigenvalues of $\left(sI-A\right)^{-1}$ associated to the
eigenvectors $v_i$.
Then we can write:
\begin{equation}
\Lambda= \left[ \begin{array}{cccccc}
\dfrac{1}{s-\lambda_{-N}}& 0 & \cdots & 0& \cdots & 0\\
0 & \dfrac{1}{s-\lambda_{-N+1}}& \cdots & \vdots& \vdots & \vdots\\
\vdots & \vdots  & \ddots & \vdots& \vdots & 0\\
\vdots & \vdots  & \cdots & \dfrac{1}{s-\lambda_0}& \vdots & \vdots\\
\vdots & \vdots  & \cdots & \vdots& \ddots & 0\\
\vdots & \vdots  & \cdots & 0& \vdots & \dfrac{1}{s-\lambda_N}
\end{array}
\right]
\end{equation}
Also define the similarity matrix $T$, considering the $2N+1$ independent eigenvectors $v_i$:
\begin{equation}
T= \left[ \begin{array}{ccccccc}
v_{-N}& \cdots  &v_{-1}&v_0&v_{1}& \cdots  &v_{N}
\end{array}
\right]
\end{equation}
And we can write:
\begin{equation}
\left(SI-A\right)^{-1}= T \Lambda T^{-1}.
\end{equation}
If we decompose $B_1$ into  directions $v_i$, i. e.,
\bean
\label{67}
B_1=\vartheta_{-N}v_{-N}+\dots+\vartheta_{-1}v_{-1}+\vartheta_0v_0+\vartheta_1v_1+\dots+\vartheta_Nv_N=T\left[ \begin{array}{c}
\vartheta_{-N} \\ \vdots \\ \vartheta_N
\end{array}\right]
\eean
then we can express the transfer function $G_{11}(s)$
as
\bean
\label{68}
\ba{rcl}
G_{11}(s)&=&C_1\left(sI-A\right)^{-1}B_1=C_1 T \Lambda T^{-1} T\left[ \begin{array}{c}
\vartheta_{-N} \\ \vdots \\ \vartheta_N
\end{array}\right]=C_1T \Lambda \left[ \begin{array}{c}
\vartheta_{-N} \\ \vdots \\ \vartheta_N
\end{array}\right]=\\
& =&C_1\left(\dfrac{\vartheta_{-N}}{s-\lambda_{-N}}v_{-N}+\dots
\dfrac{\vartheta_{-1}}{s-\lambda_{-1}}v_{-1}+\dfrac{\vartheta_0}{s}v_{0}+\dfrac{\vartheta_1}{s-\lambda_{1}}v_{1}+
\dots+\dfrac{\vartheta_{N}}{s-\lambda_{N}}v_{N}\right).
\ea
\eean
After multiplying \eqref{68} by $s$ one obtains:
\begin{equation*}
\label{70}
\ba{rcl}
K_G&=&C_1\displaystyle \lim_{s\rightarrow 0} \left(\dfrac{\vartheta_{-N}s}{s-\lambda_{-N}}v_{-N}+\dots+
\dfrac{\vartheta_{-1}s}{s-\lambda_{-1}}v_{-1}+\dfrac{\vartheta_0s}{s}v_{0}+\dfrac{\vartheta_1s}{s-\lambda_{1}}v_{1}+
\dots+\dfrac{\vartheta_{N}s}{s-\lambda_{N}}v_{N}\right)\\
&=&\vartheta_0C_1v_0=\vartheta_0,
\ea
\end{equation*}

\hfill$\Box$

because $C_1=e_1^T$ and $v_0=e_1+e_2+\dots+e_{N+1}$.
If we knew all the eigenvectors we could straightforwardly determine $\vartheta_0.$  But,  only $v_0$ is known and it is not
so immediate to compute $\vartheta_0$. The next lemma is of good help to solve this problem.
\begin{lemma}:
\label{L2}
If $A\in \cR^{n \times n}$ is a singular matrix,  $v_0$ its eigenvector associated to the zero eigenvalue and $A^Tv_0=0,$
then $v_0$ is orthogonal
to the remaining eigenvectors $v_i$, $i\neq 0,$ of $A$.
\end{lemma}
\textbf{Proof:}
Let $v_i$ with $i\neq 0$ the eigenvector of $A$ associated to the eigenvalue $\lambda_i\neq 0$. By the eigenvector
definition
\bean
\label{71}
Av_i=\lambda_iv_i\Rightarrow v_i=\dfrac{1}{\lambda_i}Av_i.
\eean
Now, using this equation, we compute the internal product between $v_i$ and $v_0$,
\bean
\label{72}
\left< v_i, v_0 \right>=v_i^Tv_0=\left(\dfrac{1}{\lambda_i}Av_i\right)^Tv_0=
\dfrac{1}{\lambda_i}v_i^TA^Tv_0=0
\eean
from which  conclude that $v_0$ and $v_i$ are orthogonal.

\hfill $ \Box $

\begin{corollary}
Consider $A$ as defined in \eqref{12}--\eqref{15}.
Its eigenvectors $v_i$, $i=-N,\ldots,-1,1,\dots,N$, are orthogonal to $v_0$.
\end{corollary}
\textbf{Proof:}
  Recalling the definitions of $v_0$ and  $A_{11}$ in equations \eqref{28} and  \eqref{12}, respectively, then:
\bea
A^Tv_0=\left[\ba{c|c}A_{11}^T & A_{21}^T \\\hline  A_{12}^T & A_{22}^T \ea\right]v_0=
\left[\ba{cc} 0_{(N+1) \times (N+1)} \\ A_{12}^Tv \ea\right]
\eea
Computing
\bea
A_{12}^Tv
&=&\left[
\ba{cccccc}
-\dfrac{c^2}{\mathcal{A} \Delta \ell} & 0 & 0 & \cdots & 0 & 0\\
[4mm]
\dfrac{c^2}{\mathcal{A} \Delta \ell} & -\dfrac{c^2}{\mathcal{A} \Delta \ell} & \vdots & \cdots & \vdots & \vdots \\
[4mm]
0 &\dfrac{c^2}{\mathcal{A} \Delta \ell} & -\dfrac{c^2}{\mathcal{A} \Delta \ell} & \cdots  & \vdots & \vdots \\
\vdots & \vdots & \ddots &  \ddots & \vdots & \vdots\\
\vdots & \vdots & \vdots & \cdots & \dfrac{c^2}{\mathcal{A} \Delta \ell} & -\dfrac{c^2}{\mathcal{A} \Delta \ell}\\
[4mm]
0 & 0 & 0 & \cdots & 0 & \dfrac{c^2}{\mathcal{A} \Delta \ell}
\ea
\right]^T
\left[\ba{c} 1 \\[5mm] 1 \\[5mm] 1 \\[4mm] \vdots \\[4mm] \vdots \\[4mm] 1 \\[5mm] 1 \ea \right]=
\eea
\bea
&= &
\left[
\ba{cccccc}
-\dfrac{c^2}{\mathcal{A} \Delta \ell} & \dfrac{c^2}{\mathcal{A} \Delta \ell}  & 0 & \cdots  & 0 & 0\\
[4mm]
0 &  -\dfrac{c^2}{\mathcal{A} \Delta \ell} & \dfrac{c^2}{\mathcal{A} \Delta \ell} & \cdots & \vdots & \vdots \\
[4mm]
 \vdots & 0 & -\dfrac{c^2}{\mathcal{A} \Delta \ell} & \ddots &   \vdots & \vdots \\
[4mm]
\vdots & \vdots &  \vdots & \ddots &   \ddots & \vdots \\
\vdots & \vdots & \vdots & \cdots &  \dfrac{c^2}{\mathcal{A} \Delta \ell} & 0\\
[4mm]
0 & 0 & 0 & \cdots &  -\dfrac{c^2}{\mathcal{A} \Delta \ell} & \dfrac{c^2}{\mathcal{A} \Delta \ell}
\ea
\right]
\left[\ba{c} 1 \\[6mm] 1 \\[6mm] 1 \\[6mm]  \vdots \\[6mm] 1 \\[6mm] 1 \ea \right]=
\eea
\bea
&= &
\left[
\ba{c}
\dfrac{c^2}{\mathcal{A} \Delta \ell}-\dfrac{c^2}{\mathcal{A} \Delta \ell}\\
\vdots\\
\dfrac{c^2}{\mathcal{A} \Delta \ell}-\dfrac{c^2}{\mathcal{A} \Delta \ell}
\ea\right]=0,
\eea
Then   by Lemma~\ref{L2} we have the expected result.

\hfill $\Box$

\begin{corollary}\label{cor5}
$K_{G}=\dfrac{c^{2}}{{\cal A}L}.$
\end{corollary}
\textbf{Proof:}
Decompose $B_1$ as
\bean
\label{73}
B_1=P_{v_0}B_1+P_{v_0^{\bot}}B_1
\eean
where $P_{w}$ is the orthogonal projection into $w$ operator and $w^{\bot}$ denotes the orthogonal complement
of $w$. From the orthogonality condition between $v_0$ and $v_i$, $i=-N,\ldots,-1, 1,\dots, N$,
\bean
\label{74}
P_{v_0}B_1=\vartheta_0v_0.
\eean
Given that
\bean
\label{75}
P_{v_0}B_1=(v_0^Tv_0)^{-1}v_0^TB_1v_0,
\eean
then
\bean
\label{76}
\vartheta_0=(v_0^Tv_0)^{-1}v_0^TB_1=
\dfrac{c^2}{\mathcal{A} (N+1)\Delta \ell}.
\eean
Replacing $\Delta \ell$ by $\dfrac{L}{N}$ we find
\bean
\label{77}
\vartheta_0=\dfrac{c^2N}{\mathcal{A} (N+1)L}.
\eean
When $N\rightarrow\infty \Rightarrow\dfrac{N}{N+1}\rightarrow 1$,
\bean
\label{78}
\vartheta_0=\dfrac{c^2}{\mathcal{A} L}.
\eean
and $K_G=\vartheta_0=\dfrac{c^2}{\mathcal{A} L}.$

\hfill$\Box$

\begin{corollary}
\label{Cor6}
Transfer function $G_{11}$ has the following form, according to Corollary~\ref{Cor3}:
\begin{equation}\label{TFG11c}
G_{11}=\dfrac{c^{2}}{\mathcal{A} L} \dfrac{ \prod_{k=1}^{\infty} \left( \dfrac{s^{2}}{z_{k}z_{-k}}+s\left( \dfrac{1}{z_{k}}+\dfrac{1}{z_{-k}}\right)+1 \right)}{s \prod_{k=1}^{\infty} \left( \dfrac{s^{2}}{\lambda_{k}\lambda_{-k}}+s \left(\dfrac{1}{\lambda_{k}}+\dfrac{1}{\lambda_{-k}}\right)+1 \right) }
\end{equation}
where $z_{k}$ are defined in Corollary~\ref{Cor2} and  $\lambda_{k}$ are defined in Corollary~\ref{Cor1}.
\end{corollary}
%=====================
\subsection{Transfer function $G_{22}$}
Next, we determine the transfer function
\begin{equation}\label{TFG22}
G_{22}=\dfrac{Y_2(s)}{U_2(s)}=\dfrac{P_{N}(s)}{Q_N(s)}
\end{equation}
when $U_{1}(s)=0.$  That is the ratio between the pressure and the massflow at the offtake node.
From \eqref{TFcomponents}:
\begin{equation}\label{82}
G_{22}(s)=C_2\left(sI-A\right)^{-1}B_2,
\end{equation}
where $C_2$ is the second column  of  $C$ and $B_2$ is the second  row of $B$.

Similarly to what happens with $G_{11}(s),$  $G_{22}(s)$ is a rational function whose poles are the eigenvalues of $A$. Following the same methodology as for  $G_{11}(s),$ we would like to calculate the zeros of $G_{22}(s)$ in order to investigate pole-zero cancelations.
\begin{theorem}\label{Th4}
The zeros of $G_{22} (s)$ are
\begin{eqnarray}
z_{ k}&=&
-\dfrac{f_c c^2 Q_{m}}{4 \mathcal{D} \mathcal{A} P_{m}} \pm j
\sqrt{\left(2\dfrac{c}{\Delta \ell}\sin\left(\dfrac{(2k-1)\pi}{2N+1}\right)\right)^2
-\left(\dfrac{f_c c^2 Q_{m}}{4 \mathcal{D} \mathcal{A} P_{m}}\right)^2},\quad k= 1,\dots, N\label{87}
\end{eqnarray}
\end{theorem}
\textbf{Proof:}
The proof is very similar to the one of Theorem~\ref{Th2}.
Again, we recall the result from \cite[pp. 284]{Carvalho:1993} that states that the zeros of the transfer function \eqref{TFG22} are the zeros of the following polynomial:
\begin{equation}
\left|\begin{array}{c|c}
sI_{(2N+1)} -A & -B_2 \\\hline
C_2 & 0
\end{array}
\right|=0
\end{equation}
Recall the definitions of  $C_2=e_{N+1} \in \cR^{(2N+1)}$ and $B_2=-\dfrac{c^2}{{\cal A} \Delta \ell }e_{N+1}\cR^{(2N+1)}$ and the proof follows exactly as for Theorem~\ref{Th2}.

We have:
\begin{eqnarray}
\left|\begin{array}{c|c}
sI_{(2N+1)} -A & e_{N+1}\\\hline
e_{N+1}^T& 0
\end{array}
\right|&=& 0 \quad \Leftrightarrow\\[2mm]
\Leftrightarrow \left|\begin{array}{c|c}
\left[ \begin{array}{c|c} sI_{N+1}-A_{11} & -A_{12} \\\hline -A_{21} & sI_N-A_{22}\end{array}\right] & \begin{array}{c} 0 \\  \vdots \\ 0 \\ \hline 1 \\ \vdots \\0 \end{array}\\[2mm]\hline
\begin{array}{cccc|cccc}0 & \cdots &  \cdots& 0 & 1 & \cdots & \cdots & 0\end{array}&  0
\end{array}
\right|
\end{eqnarray}
From the definition of $A_{11}$ and $A_{22}$ defined in \eqref{12} and \eqref{15}
\begin{eqnarray}
&=& 0\quad \Leftrightarrow\\[2mm]
\Leftrightarrow\left|\begin{array}{c|c}
\left[ \begin{array}{c|c} sI_{N+1} & -A_{12} \\ \hline -A_{21} & (s+2\alpha)I_N\end{array}\right] & \begin{array}{c} 0 \\  \vdots \\ 0 \\ \hline 1 \\ \vdots \\0 \end{array}\\[2mm]\hline
\begin{array}{ccc|ccc}0 &  \cdots& 0 & 1 & \cdots & 0\end{array}&  0
\end{array}
\right|&=& 0.
\end{eqnarray}
We develop this determinant first along the last column and next along the last row, and obtain:
\begin{eqnarray}
\left| \begin{array}{c|c} sI_{N} & -\bar{A}_{12} \\[4mm] \hline \\ -\bar{A}_{21} & (s+2\alpha)I_N\end{array}
\right|&=& 0
\end{eqnarray}
where $\bar{A}_{12}$ denotes matrix $A_{12}$ without the $1$st column and  $\bar{A}_{21}$ denotes matrix $A_{21}$ without the $1$st row.
Next, as $sI_{N}\left(-\bar{A}_{12} \right)=\left(-\bar{A}_{12} \right)sI_{N} $, we apply again
Fact 2.13.10 in \cite[pp. 62--63]{Bernstein}, that states that for the arbitrary matrices
$\mathbf{A},\;\mathbf{B},\;\mathbf{C}$ and $\mathbf{D} \in \cR^{N \times N}$ such that
$\mathbf{A}\mathbf{B}=\mathbf{B}\mathbf{A}$ then
\bea
\det\left[\ba{c|c} \mathbf{A} & \mathbf{B} \\\hline \mathbf{C} & \mathbf{D} \ea \right]=
\det(\mathbf{D}\mathbf{A}-\mathbf{C}\mathbf{B})
\eea
and
\begin{eqnarray}
\left| \begin{array}{c|c} sI_{N} & -\bar{A}_{12} \\ \hline  -\bar{A}_{21} & (s+2\alpha)I_N\end{array}
\right|&=& \left| s(s+2\alpha)I_N-\bar{A}_{21}\bar{A}_{12}\right|\\[2mm]
&=& \left|\left(s^2+2\alpha s+\alpha^2 \right) I_N-\left( \bar{A}_{21}\bar{A}_{12} +\alpha^2\right) \right|\\[2mm]
&=&\left|\left(s+\alpha\right)^2  I_N-\left( \bar{A}_{21}\bar{A}_{12} +\alpha^2\right) \right|.
\end{eqnarray}
We do the usual change of variable
\begin{equation}\label{chvar2}
\mathcal{S}=(s+\alpha)^2
\end{equation}
and then can write:
\begin{eqnarray}
&=&\left|  \mathcal{S}I_N-\left( \bar{A}_{21}\bar{A}_{12} +\alpha^2\right) \right|
\end{eqnarray}
Next, we calculate  the spectrum of matrix $\left( \bar{A}_{21}\bar{A}_{12} +\alpha^2\right),$ that is:
\[
\Lambda \left( \bar{A}_{21}\bar{A}_{12} +\alpha^2\right)=\Lambda \left( \bar{A}_{21}\bar{A}_{12} \right)+\alpha^2
\]
Now, we calculate the product:
\begin{eqnarray*}
\bar{A}_{21}\bar{A}_{12}& =&  \left[
\ba{ccccc}
- \dfrac{\mathcal{A}}{\Delta \ell} & 0 &\cdots & 0 & 0 \\
[4mm]
\dfrac{\mathcal{A}}{\Delta \ell} & - \dfrac{\mathcal{A}}{\Delta \ell} & \cdots & \vdots & \vdots \\
0 & \ddots & \ddots & \vdots & \vdots\\
0 & 0 & \cdots & \dfrac{\mathcal{A}}{\Delta \ell} & - \dfrac{\mathcal{A}}{\Delta \ell}\\
0 & 0 & \cdots & 0&\dfrac{\mathcal{A}}{\Delta \ell} \\
\ea
\right]\left[
\ba{ccccccc}
\dfrac{c^2}{\mathcal{A} \Delta \ell} & -\dfrac{c^2}{\mathcal{A} \Delta \ell} & 0 & 0 & \cdots & 0 & 0\\
[4mm]
0 &\dfrac{c^2}{\mathcal{A} \Delta \ell} & -\dfrac{c^2}{\mathcal{A} \Delta \ell} & 0 & \cdots & \vdots & \vdots \\
[4mm]
\vdots & \vdots & \vdots & \ddots &  \cdots & \vdots & \vdots\\
\vdots &\vdots & \vdots & \ddots & \cdots & -\dfrac{c^2}{\mathcal{A} \Delta \ell} &0\\
0 &0&0&0 &\cdots &\dfrac{c^2}{\mathcal{A} \Delta \ell} & -\dfrac{c^2}{\mathcal{A} \Delta \ell}
\ea
\right]\\[4mm]
& =& (\dfrac{c}{\Delta \ell})^2  \left[\ba{cccccc} -2 & 1 & 0 & \cdots & 0 & 0 \\ 1 & -2 & 1 & \cdots & 0 & 0 \\ 0 & 1 & -2 & \ddots & 0 & 0 \\
\vdots & \vdots & \ddots & \ddots & \ddots & \vdots \\ 0 & 0 & 0 & \ddots & -2 & 1 \\ 0 & 0 & 0 & \cdots & 1 & -1\ea\right]\in \cR^N
\end{eqnarray*}
Then
\[ \left|sI_N- \bar{A}_{21}\bar{A}_{12}\right|= (\dfrac{c}{\Delta \ell})^2
\Lambda\left(M_N\right)
\]
where
\begin{equation}
\label{61b}
M_N=\left[\ba{cccccc} -2 & 1 & 0 & \cdots & 0 & 0 \\ 1 & -2 & 1 & \cdots & 0 & 0 \\ 0 & 1 & -2 & \ddots & 0 & 0 \\
\vdots & \vdots & \ddots & \ddots & \ddots & \vdots \\ 0 & 0 & 0 & \ddots & -2 & 1 \\ 0 & 0 & 0 & \cdots & 1 & -1\ea\right]
\in \cR^{N \times N}.
\end{equation}
From \cite[pp. 72]{Yueh:2005} \begin{equation}
\label{62b}
\Lambda\left(M_N\right)=-2+2\cos\left(\dfrac{\left( 2k-1\right) \pi}{2N+1}\right),\quad k=1,2,3,\ldots,N.
\end{equation}
and from here the proof follows exactly as for the calculus of the zeros of $G_{11}(s),$ and we can see that are the same.

Likewise follows for the asymptotic case.
As we can see from the definition of the transfer function \eqref{43} and \eqref{82} as well as from the definition of the $B_{i},C_{i},i=1,2$ we have
\bean
\label{80}
G_{11}(s)=-G_{22}(s)
\eean
Therefore, its zeros will be necessarily coincident.

$\hfill\Box$

%============================================
\subsection{Transfer function $G_{12}$}
Acccording to \eqref{TFcomponents}, consider now the transfer function
\begin{equation}\label{TFG12}
G_{12}(s)=\dfrac{Y_1(s)}{U_2(s)}=\dfrac{P_{0}(s)}{Q_{N+1}(s)}
\end{equation}
with $U_1(s)=Q_0(s)=0,$
or equivalently:
\bean
\label{107}
G_{12}(s)&=&C_1\left(sI-A\right)^{-1}B_2,
\eean
where $C_1$ is the first column  of  $C$ and $B_2$ is the second  row of $B$.  \\
\begin{theorem}\label{Th5}
The transfer function   $G_{12}(s)$ has no zeros.
\end{theorem}
\textbf{Proof:}
According to \cite[pp. 284]{Carvalho:1993}, we have that the zeros of the transfer function \eqref{TFG12} are the zeros of the following polynomial:
\begin{equation*}
\left|\begin{array}{c|c}
sI_{(2N+1)} -A & -B_2 \\\hline
C_1 & 0
\end{array}
\right|
\end{equation*}
and recalling that $C_1=e_1 \in \cR^{(2N+1)}$ and $B_2=-e_{N+1}\cR^{(2N+1)}$, we have:
\begin{eqnarray}
\left|\begin{array}{c|c}
sI_{(2N+1)} -A & e_{N+1} \\\hline
e_1^T& 0
\end{array}
\right|&=& 0 \quad \Leftrightarrow \nonumber \\[2mm]
\Leftrightarrow \left|\begin{array}{c|c}
\left[ \begin{array}{c|c} sI_{N+1}-A_{11} & -A_{12} \\\hline -A_{21} & sI_N-A_{22}\end{array}\right] & \begin{array}{c} 0 \\ \vdots \\ 0 \\ 1 \\ 0 \\ \vdots \\ 0 \end{array}\\[2mm]\hline
\begin{array}{cccccccc}1 & 0 &  \cdots & \cdots & \cdots & \cdots & 0\end{array}&  0
\end{array}
\right|
\end{eqnarray}
And from de definition of $A_{11}$ and $A_{22}$ in \eqref{12} and \eqref{15}
\begin{eqnarray}
&=& 0\quad \Leftrightarrow \nonumber \\[2mm]
\Leftrightarrow\left|\begin{array}{c|c}
\left[ \begin{array}{c|c} sI_{N+1} & -A_{12} \\\hline -A_{21} & (s-\alpha)I_N\end{array}\right] & \begin{array}{c} 0 \\ \vdots \\ 0 \\ 1 \\ 0 \\ \vdots \\ 0 \end{array}\\[2mm]\hline
\begin{array}{cccccc}1 & 0 & \cdots& \cdots & \cdots  & 0\end{array}&  0
\end{array}
\right|&=& 0 \label{mz(N+1)}
\end{eqnarray}
We develop this determinant first along the last column and next along the last row, and obtain:
\begin{equation}
\label{mxNiter1}
\begin{array}{rcl}
\small
\left| \begin{array}{ccccccc|cccccccc} 0 & 0 & 0 & \cdots & 0 &  0 &0 & \gamma & 0 & 0 & \cdots & 0 & 0 & 0   \\
s & 0 & 0 & \cdots & 0& 0 & 0 & -\gamma & \gamma &0 & \cdots & \vdots & \vdots &\vdots   \\
0 & s & 0 & \cdots& 0 & 0 & 0 & 0 & -\gamma & \gamma &  \cdots & \vdots & \vdots &\vdots   \\
\vdots & \vdots & \vdots & & \vdots & \vdots & \vdots &  \vdots  &\vdots  & \ddots  &  \ddots & \vdots & \vdots & \vdots    \\
0 & 0 & 0 & \ddots&s & 0 & 0 &  0 &0 & 0 &    \cdots & -\gamma & \gamma &0   \\
0 & 0 & 0 & \cdots& 0& s & 0 &  0 &0 & 0 &    \cdots & 0 & -\gamma & \gamma    \\
\hline
\varrho  & 0 & 0 & \cdots& 0 & 0 & 0 &(s+2\alpha) &0 & 0 &  \cdots & 0 & 0 & 0   \\
-\varrho  &  \varrho & 0 & \cdots& \vdots & \vdots  & \vdots  &  0 &(s+2\alpha)  & 0 &   \cdots & \vdots & \vdots & \vdots   \\
0 & -\varrho  &  \varrho &  \cdots& \vdots  & \vdots  &  \vdots  &0 &0 & (s+2\alpha)  &  \cdots & \vdots & \vdots & \vdots  \\
\vdots & \vdots & \vdots & \ddots & \ddots & \vdots &  \vdots  & \vdots &   \vdots  &  \vdots & \ddots & \vdots & \vdots & \vdots    \\
\vdots  & \vdots    & \vdots  &  \cdots& -\varrho &  \varrho & 0        &\vdots  & \vdots   &\vdots  &   \cdots  &  0 & (s+2\alpha) & 0   \\
0 & 0  & 0 &  \cdots&        0 & -\varrho &  \varrho & 0 &0 & 0 & \cdots &   0  &  0 & (s+2\alpha) 
\end{array}
\right|&=& 0
\normalsize
\end{array}
\end{equation}
$\gamma= -\dfrac{c^2}{{\cal A}\Delta \ell}$ and $\varrho = -\dfrac{{\cal A}}{\Delta \ell}.$\\

We don't worry about the signs of the cofactor, since our aim is to determine the zeros of the determinant.  Now, we develop this determinant first along  the first line, and we obtain:
\begin{eqnarray*}
\label{mxNiter2}
\small
\gamma \left| \begin{array}{ccccccc|cccccc}
s & 0 & 0 & \cdots & 0& 0 & 0  & \gamma &0 & \cdots & 0 & 0 &0   \\
0 & s & 0 & \cdots&  \vdots &  \vdots &  \vdots &  -\gamma & \gamma &  \cdots &  \vdots &  \vdots & \vdots   \\
\vdots & \vdots & \ddots & \cdots &  \vdots & \vdots &  \vdots &  \vdots  & \ddots  &  \ddots & \vdots & \vdots & \vdots    \\
 \vdots &  \vdots &  \vdots & \cdots&s & 0 & 0 &   \vdots  &  \vdots &    \cdots & -\gamma & \gamma &0   \\
0 & 0 & 0 & \cdots& 0& s & 0 &  0 &0  &    \cdots & 0 & -\gamma & \gamma    \\
\hline
\varrho  & 0 & 0 & \cdots& 0 &  0 & 0 &   0 & 0 &  \cdots & 0 & 0 & 0   \\
-\varrho  &  \varrho & 0 & \cdots&  \vdots&  \vdots & 0 &  (s+2\alpha)  & 0 &   \cdots &  \vdots &  \vdots &  \vdots   \\
0 & -\varrho  &  \varrho &  \cdots& \vdots  &  \vdots &  0 &0  & (s+2\alpha)  &  \cdots &  \vdots &  \vdots &  \vdots   \\
\vdots & \vdots & \vdots & \ddots & \ddots & \vdots & \vdots  & \vdots & \vdots  &  \cdots & \ddots & \vdots & \vdots    \\
 \vdots &  \vdots  &  \vdots &  \cdots& -\varrho &  \varrho & 0        & \vdots &  \vdots  &   \cdots  &   \vdots & (s+2\alpha) & 0   \\
0 & 0  & 0 &  \cdots&        0 & -\varrho &  \varrho & 0 &0  &\cdots  &  0   &  0 & (s+2\alpha) 
\end{array}
\right|&=& 0
\normalsize
\end{eqnarray*}
Next, we develop along column--$N$:
\begin{eqnarray*}
\label{mx(N-1)iter1}
\small
\gamma \varrho \left| \begin{array}{cccccc|ccccccc}
s & 0 & 0 & \cdots & 0& 0 &  \gamma &0 & \cdots & 0 & 0 &0   \\
0 & s & 0 & \cdots&  \vdots &  \vdots  &  -\gamma & \gamma &  \cdots &  \vdots &  \vdots & \vdots   \\
\vdots & \vdots & \vdots & \ddots & \vdots  & \vdots & \vdots & \vdots & \ddots & \ddots & \vdots & \vdots    \\
 \vdots &  \vdots &  \vdots & \cdots&s & 0  &   \vdots  &  \vdots &    \cdots & -\gamma & \gamma &0   \\
0 & 0 & 0 & \cdots& 0& s  &  0 &0  &    \cdots & 0 & -\gamma & \gamma    \\
\hline
\varrho  & 0 & 0 & \cdots& 0  & 0 &   0 & 0 &  \cdots & 0 & 0 & 0   \\
-\varrho  &  \varrho & 0 & \cdots&  \vdots&   \vdots &  (s+2\alpha)  & 0 &   \cdots &  \vdots &  \vdots &  \vdots   \\
0 & -\varrho  &  \varrho &  \cdots& 0  &  0 &0  & (s+2\alpha)  &  \cdots &  \vdots &  \vdots &  \vdots   \\
\vdots & \vdots & \vdots & \ddots & \ddots & \vdots & \vdots & \vdots &  \cdots & \ddots & \vdots & \vdots    \\
0 & 0  & 0 &  \cdots& -\varrho &  \varrho         &0 & 0  &   \cdots  &  0 & (s+2\alpha) & 0
\end{array}
\right|&=& 0
\normalsize
\end{eqnarray*}
Here the matrix is of dimension $2(N-1)\times 2(N-1).$

Swap the first row of blocks with the second one:
\begin{eqnarray*}
\label{mx(N-1)iter2}
\small
\gamma \varrho \left| \begin{array}{cccccc|ccccccc}
\varrho  & 0 & 0 & \cdots& 0  & 0 &   0 & 0 &  \cdots & 0 & 0 & 0   \\
-\varrho  &  \varrho & 0 & \cdots&  \vdots&   \vdots &  (s+2\alpha)  & 0 &   \cdots &  \vdots &  \vdots &  \vdots   \\
0 & -\varrho  &  \varrho &  \cdots&  \vdots  &   \vdots &0  & (s+2\alpha)  &  \cdots &  \vdots &  \vdots &  \vdots   \\
\vdots & \vdots & \vdots & \ddots & \ddots & \vdots & \vdots & \vdots &  \cdots & \ddots & \vdots & \vdots    \\
0 & 0  & 0 &  \cdots& -\varrho &  \varrho         &0 & 0  &   \cdots  &  0 & (s+2\alpha) & 0   \\
\hline
s & 0 & 0 & \cdots & 0& 0 &  \gamma &0 & \cdots & 0 & 0 &0   \\
0 & s & 0 & \cdots&  \vdots &  \vdots  &  -\gamma & \gamma &  \cdots &  \vdots &  \vdots & \vdots   \\
\vdots & \vdots & \vdots & \ddots & \vdots  & \vdots & \vdots & \vdots & \ddots & \ddots & \vdots & \vdots    \\
 \vdots &  \vdots &  \vdots & \cdots&s & 0  &   \vdots  &  \vdots &    \cdots & -\gamma & \gamma &0   \\
0 & 0 & 0 & \cdots& 0& s  &  0 &0  &    \cdots & 0 & -\gamma & \gamma
\end{array}
\right|&=& 0
\normalsize
\end{eqnarray*}
 Swap the first column of blocks with the second one:
\begin{eqnarray*}
\small
\gamma \varrho \left| \begin{array}{cccccc|ccccccc}
  0 & 0 &  \cdots & 0 & 0 & 0 &\varrho  & 0 & 0 & \cdots& 0  & 0   \\
 (s+2\alpha)  & 0 &   \cdots &  \vdots &  \vdots &  \vdots  & -\varrho  &  \varrho & 0 & \cdots&  \vdots&   \vdots \\
0  & (s+2\alpha)  &  \cdots &  \vdots &  \vdots &  \vdots  &0 & -\varrho  &  \varrho &  \cdots&  \vdots  &   \vdots  \\
\vdots & \vdots & \cdots & \ddots & \vdots & \vdots  & \vdots & \vdots & \vdots & \ddots & \ddots  & \vdots    \\
0 & 0  &   \cdots  &  0 & (s+2\alpha) & 0  & 0 & 0  & 0 &  \cdots& -\varrho &  \varrho  \\
\hline
\gamma &0 & \cdots & 0 & 0 &0 & s & 0 & 0 & \cdots & 0& 0   \\
  -\gamma & \gamma &  \cdots &  \vdots &  \vdots & \vdots  & 0 & s & 0 & \cdots&  \vdots &  \vdots   \\
  \vdots & \vdots  &  \ddots & \ddots & \vdots & \vdots & \vdots & \vdots & \vdots & \ddots & \vdots & \vdots \\
  \vdots   &  \vdots &    \cdots & -\gamma & \gamma &0   &  \vdots &  \vdots &  \vdots & \cdots&s & 0  \\
  0 &0  &    \cdots & 0 & -\gamma & \gamma    &0 & 0 & 0 & \cdots& 0& s
\end{array}
\right|&=& 0
\normalsize
\end{eqnarray*}
Again develop the determinant along the first row:
\begin{eqnarray*}
\small
\gamma \varrho^2 \left| \begin{array}{cccccc|ccccc}

 (s+2\alpha)  & 0 &   \cdots & 0 & 0 & 0  &   \varrho & 0 & \cdots& 0&  0 \\
0  & (s+2\alpha)  &  \cdots & \vdots  &  \vdots &  \vdots   & -\varrho  &  \varrho &  \cdots&  \vdots  &   \vdots  \\
\vdots & \vdots & \vdots & \ddots & \vdots & \vdots  & \vdots   & \vdots & \ddots& \ddots & \vdots \\
0 & 0  &   \cdots  &  0 & (s+2\alpha) & 0  & 0  & 0 &  \cdots& -\varrho &  \varrho   \\
\hline
\gamma &0 & \cdots & 0 & 0 &0 &  0 & 0 & \cdots & 0& 0   \\
  -\gamma & \gamma &  \cdots &  \vdots &  \vdots & \vdots  &  s & 0 & \cdots&  \vdots &  \vdots   \\
\vdots & \vdots & \ddots & \ddots & \vdots & \vdots & \vdots & \vdots & \ddots & \vdots & \vdots  \\
   \vdots  &  \vdots &    \cdots & -\gamma & \gamma &0    &  \vdots &  \vdots & \cdots&s & 0  \\
  0 &0  &    \cdots & 0 & -\gamma & \gamma   & 0 & 0 & \cdots& 0& s
\end{array}
\right|&=& 0
\normalsize
\end{eqnarray*}
Again along column--$(N-1)$
\begin{eqnarray*}
\small
\gamma^2  \varrho^2 \left| \begin{array}{ccccc|ccccc}
 (s+2\alpha)  & 0 &   \cdots & 0 & 0  &   \varrho & 0 & \cdots& 0&  0 \\
0  & (s+2\alpha)  &  \cdots &  \vdots &  \vdots    & -\varrho  &  \varrho &  \cdots&  \vdots  &   \vdots  \\
\vdots & \vdots & \cdots & \ddots & \vdots & \vdots & \vdots & \ddots & \ddots & \vdots \\
0 & 0  &   \cdots  &  0 & (s+2\alpha)   & 0  & 0 &  \cdots& -\varrho &  \varrho   \\
\hline
\gamma &0 & \cdots & 0 & 0  &  0 & 0 & \cdots & 0& 0   \\
  -\gamma & \gamma &  \cdots &  \vdots &  \vdots   &  s & 0 & \cdots&  \vdots &  \vdots   \\
\vdots & \vdots & \ddots & \ddots & \vdots & \vdots & \vdots & \ddots & \vdots & \vdots \\
  0  & 0 &    \cdots & -\gamma  & \gamma    & 0 & 0 & \cdots&s & 0
\end{array}
\right|&=& 0
\normalsize
\end{eqnarray*}
Therefore, we find a pattern.  To write the pattern, we define:
\begin{eqnarray*}
\tilde{I}_N &=& \mbox{identity matrix whose first row and column--$N$ are all zeros}\\
\bar{A}_{12}&=& \mbox{matrix $A_{12}$ without the last row}\\
\bar{A}_{{12}_N}&=& \mbox{ matrix of order $N$ and with the same pattern as $\bar{A}_{12}$}\\
\bar{A}_{21}&=& \mbox{ matrix $A_{21}$ without the first column}\\
\bar{A}_{{21}_N}&=& \mbox{ matrix of order $N$ and with the same pattern as $\bar{A}_{21}$}\\
\end{eqnarray*}
with this notation, we can write the determinant \eqref{mxNiter1} as:
\begin{eqnarray*}
\mbox{with }i=0,  \mbox{it becomes} \ \left| \begin{array}{c|c}
\\
s \tilde{I}_{N-i} & \bar{A}_{{12}_{N-i}}\\[2mm]
\hline\\
\bar{A}_{{21}_{N-i}} & (s+2\alpha) I_{N-i}
\\[2mm]
\end{array}\right|
\end{eqnarray*}
Define an iteration as:
\begin{enumerate}
\item Develop the determinant in cofactors along the first row
\item Develop the determinant in cofactors along column $(N-i)$
\item Switch the first row of blocks with the second one
\item Switch the first column of blocks with the second one
\end{enumerate}
Then, we obtain:
\begin{eqnarray*}
\gamma^{i+1}\varrho^{i+1} \left| \begin{array}{c|c}
\\
(s+2\alpha)  \tilde{I}_{N-(i+1)} & \bar{A}_{{21}_{N-(i+1)}}\\[2mm]
\hline\\
\bar{A}_{{12}_{N-(i+1)}} & sI_{N-(i+1)}
\\[2mm]
\end{array}\right|
\end{eqnarray*}
Iterate again and obtain:
\begin{eqnarray*}
\gamma^{i+2}\varrho^{i+2} \left| \begin{array}{c|c}
\\
s \tilde{I}_{N-(i+2)} & \bar{A}_{{12}_{N-(i+2)}}\\[2mm]
\hline\\
\bar{A}_{{21}_{N-(i+2)}} & (s+2\alpha) I_{N-(i+2)}
\\[2mm]
\end{array}\right|
\end{eqnarray*}
Also, considering $N=1$ we write \eqref{mxNiter1} as:
\begin{eqnarray}
\left| \begin{array}{cc|c|c}
s & 0 & -\gamma & 0 \\
0 & s & \gamma & 1 \\\hline
\varrho & -\varrho & (s-\alpha ) & 0 \\\hline
1 & 0 & 0 & 0
\end{array}\right|
\end{eqnarray}
Similarly to what we have done for the general case, we develop the determinant first along the last column
\begin{eqnarray*}
\left| \begin{array}{cc|c}
s & 0 & -\gamma  \\ \hline
\varrho & \varrho & (s-\alpha )  \\ \hline
1 & 0 & 0
\end{array}\right|
\end{eqnarray*}
and next along the last row:
\begin{eqnarray*}
\left| \begin{array}{c|c}
0 & -\gamma  \\ \hline
\varrho & (s-\alpha )
\end{array}\right|&=&-\varrho\gamma \neq 0,
\end{eqnarray*}
and the proof that $G_{12}(s)$ has no zeros is complete.
\hfill $\Box$

$G_{12}(s)$ is a rational function
whose poles are all the eigenvalues of $A$, since this transfer functions has no zeros. Therefore, we can write:
\begin{corollary}\label{Cor7}
$G_{12}(s)$ is given by
\begin{equation}
\label{d2}
G_{12}(s)=\dfrac{K_G}{s  \left( \dfrac{s^{2}}{\lambda_{k}\lambda_{-k}}+s \left(\dfrac{1}{\lambda_{k}}+\dfrac{1}{\lambda_{-k}}\right)+1 \right) }
\end{equation}
where $\lambda_{\pm k}$ is given by \eqref{assymeigen}.
\end{corollary}
%========================
%========================
\section{Approximated Transfer functions}\label{}
In this section we propose some approximations for  the models of  the transfer functions.

From Corollaries~\ref{Cor1}--\ref{Cor2} and Theorem \ref{Th3}, $G_{ij}(s)$, $i,\;j=1,\;2,$ are meromorphic functions given by

\begin{equation}
 \begin{array}{l}\label{29}
     \hspace{-5mm}
     G_{11}(s) = \dfrac{K_G}{s}\displaystyle\prod_{k=1}^{\infty}K_k
     \dfrac{\left(s+\alpha\right)^2+(2k-1)^2\omega_0^2-\alpha^2}{\left(s+\alpha\right)^2+4k^2\omega_0^2-\alpha^2} \\
     \hspace{-5mm}
     G_{21}(s) = \dfrac{K_G}{s}\displaystyle\prod_{k=1}^{\infty}\dfrac{4k^2\omega_0^2}{\left(s+\alpha\right)^2+4k^2\omega_0^2-\alpha^2}\\
     \hspace{-5mm}
     G_{22}(s) = -G_{11}(s)\\
      \hspace{-5mm}
    G_{12}(s) = -G_{21}(s)
 \end{array}
\end{equation}
where
\begin{equation}\label{30}
  \omega_0=\dfrac{\pi}{2T_d}
\end{equation}
and
\begin{equation}\label{31}
  \begin{array}{rcl}
     K_k & = & \left(\dfrac{2k}{2k-1}\right)^2
  \end{array}
\end{equation}
Natural gas is highly pressurized in transportation networks in order to expedite its flow. To ensure this, it must compressed  periodically along the pipe. This is accomplished by compressor stations, which are usually placed at 60 Km to 250 Km intervals along the pipeline. As a result, the frequency $\omega_0$ always remains much greater than $\alpha.$  Taking this into account as well as the requirement that its factors have a DC gain set to 1, we define $ \widehat{K}_k$,   and thence  can approximate $G_{ij}(s)$, $i,j=1,2$ by
\begin{equation}\label{32}
  \begin{array}{l}
     \widehat{G}_{11}(s)  = \dfrac{K_G}{s} \displaystyle\prod_{k=1}^{\infty}\widehat{K}_k\dfrac{(s+\alpha)^2+(2k-1)^2\omega_0^2}
     {(s+\alpha)^2+4k^2\omega_0^2}\\[5mm]
     \widehat{G}_{21}(s) =  \dfrac{K_G}{s}\displaystyle\prod_{k=1}^{\infty} \dfrac{\alpha^2+4k^2\omega_0^2}{(s+\alpha)^2+4k^2\omega_0^2}\\[5mm]
     \widehat{G}_{12}(s)=-\widehat{G}_{21}(s)\\
     \widehat{G}_{22}(s)=-\widehat{G}_{11}(s)
  \end{array}
\end{equation}
with
\begin{equation}\label{33}
  \widehat{K}_k=\dfrac{\alpha^2+4k^2\omega_0^2}{\alpha^2+(2k-1)^2\omega_0^2}.
\end{equation}
If we define $\mathrm{S}=s+\alpha$ we can write
\begin{equation}\label{34}
  \begin{array}{l}
      \hspace{-5mm}
      \widehat{G}_{11}(s)=\bar{G}_{11}(\mathrm{S})= \displaystyle \dfrac{K_G}{\mathrm{S}-\alpha}\prod_{k=1}^{\infty}\widehat{K}_k\dfrac{\mathrm{S}^2+(2k-1)^2\omega_0^2}{\mathrm{S}^2+4k^2\omega_0^2}\\[5mm]
      \hspace{-5mm}
      \widehat{G}_{21}(s)=\bar{G}_{21}(\mathrm{S})=\dfrac{K_G}{\mathrm{S}-\alpha}\displaystyle\prod_{k=1}^{\infty} \dfrac{\alpha^2+4k^2\omega_0^2}{\mathrm{S}^2+4k^2\omega_0^2}.
   \end{array}
\end{equation}

Theorem~\ref{Th5} considers auxiliary functions that lead to significant simplification in the  representation of $\bar{G}_{11}(\mathrm{S})$ and $\bar{G}_{21}(\mathrm{S}).$  However, before stating Theorem~\ref{Th5}  we need to prove some intermediate results:

\begin{theorem}\label{th2D0}
The function
\begin{equation}\label{auxf}
f(s)=\dfrac{e^{-sT_d}}{1-e^{-2sTd}}
\end{equation}
be expanded as
\begin{equation}\label{auxfexpansion}
   f(s)=\sum_{k=-\infty}^{\infty}\dfrac{a_k}{s-\lambda_k}
\end{equation}
where $a_k$ is the residual of $f(s)$ at $s=\lambda_k$, i. e.
$$  a_k=\dfrac{(-1)^k}{2T_d}.$$
\end{theorem}
Since condition \eqref{c1} holds (see appendix \ref{RationalExpansionMF}) then the expansion \eqref{auxfexpansion} exits and the residuals $a_k$ are given by
\textbf{Proof:}
\begin{eqnarray*}
   a_k&=&\dfrac{1}{2j\pi}\oint_{\mathcal{C}_k}f(s)ds=\lim_{s\rightarrow \lambda_k}(s-\lambda_k)f(s)=
   \lim_{s\rightarrow \lambda_k}\dfrac{(s-\lambda_k)e^{-Tds}}{1-e^{-2T_ds}}\\
   &=&
   \lim_{s\rightarrow \lambda_k}\dfrac{e^{-sT_d}-(s-\lambda_k)T_de^{-sT_d}}{2T_de^{-2sT_d}}=\dfrac{e^{\lambda_kT_d}}{2T_d}=\dfrac{e^{\frac{kj\pi}{T_d}T_d}}{2Td}=
   \dfrac{e^{jk\pi}}{2T_d}=\dfrac{(-1)^k}{2T_d}
\end{eqnarray*}
$\hfill\Box$

\begin{theorem}\label{th2D}
The function
\begin{equation}\label{auxg} \end{equation}
$$g(s)=\dfrac{\displaystyle\prod_{k=-\infty,k\neq0}^{\infty}\left(-jk\dfrac{\pi}{T_d}\right)}{\displaystyle\prod_{k=-\infty}^{\infty}(s-jk\dfrac{\pi}{T_d})}.$$
can be written as:    $$2T_d\dfrac{e^{-T_ds}}{1-e^{-2T_ds}}=2T_df(s).$$
\end{theorem}
\textbf{Proof:}
Since $g(s)$ is  proper it can be expanded as a Laurent series
\begin{equation*}
   g(s)=\sum_{k=-\infty}^{\infty}\dfrac{b_k}{s-\lambda_k}
\end{equation*}
where $b_k$ is the residual of $g(s)$ at $s=\lambda_k=jk\dfrac{\pi}{T_d}$. In order to compute this residual we rewrite $g(s)$ as
\begin{equation*}
  g(s)=\lim_{M\rightarrow\infty}
  \dfrac{\displaystyle\prod_{k=-M,k\neq0}^{M}\left(-jk\dfrac{\pi}{T_d}\right)}{\displaystyle\prod_{k=-M}^{M}(s-jk\dfrac{\pi}{T_d})}.
\end{equation*}
The residual $b_k$ is then given by
\begin{equation*}
   b_k=\dfrac{1}{2j\pi}\oint_{\mathcal{C}_k}g(s)ds=\lim_{s\rightarrow \lambda_k }(s-\lambda_k)g(s)= \lim_{M\rightarrow\infty}\dfrac{\displaystyle\prod_{m=-M,m\neq0}^{M}\left(jm\dfrac{\pi}{T_d}\right)} {\displaystyle\prod_{m=-M}^{k-1}\left(j(k-m)\dfrac{\pi}{T_d}\right)\prod_{m=k+1}^{M}\left(j(k-m)\dfrac{\pi}{T_d}\right)}
\end{equation*}
Given that
\begin{equation*}
   \left(-jm\dfrac{\pi}{T_d}\right) \left(jm\dfrac{\pi}{T_d}\right)=-j^2\left(\dfrac{m\pi}{T_d}\right)^2=\left(\dfrac{m\pi}{T_d}\right)^2
\end{equation*}
we can rewrite $b_k$ as
\begin{equation*}
   b_k=\lim_{M\rightarrow\infty}\dfrac{\displaystyle\prod_{m=1}^{M}\left(m\dfrac{\pi}{T_d}\right)^2} {\displaystyle\prod_{m=-M}^{k-1}\left(j(k-m)\dfrac{\pi}{T_d}\right)\prod_{m=k+1}^{M}\left(j(k-m)\dfrac{\pi}{T_d}\right)}.
\end{equation*}
If now we replace define $\ell=k-m$ we get
\begin{eqnarray*}
   b_k &= & \lim_{M\rightarrow\infty}\dfrac{\displaystyle\prod_{m=1}^{M}\left(m\dfrac{\pi}{T_d}\right)^2} {\displaystyle\prod_{\ell=1}^{k+M}\left(j\ell\dfrac{\pi}{T_d}\right)\prod_{\ell=-(M-k)}^{-1}\left(j\ell\dfrac{\pi}{T_d}\right)}=
   \lim_{M\rightarrow\infty}\dfrac{\displaystyle\prod_{m=1}^{M}\left(m\dfrac{\pi}{T_d}\right)^2} {\displaystyle j^{2M}\prod_{\ell=1}^{k+M}\left(\ell\dfrac{\pi}{T_d}\right)\prod_{\ell=-(M-k)}^{-1}\left(\ell\dfrac{\pi}{T_d}\right)}\\
   &=& \lim_{M\rightarrow\infty}\dfrac{\displaystyle\prod_{m=1}^{M}\left(m\dfrac{\pi}{T_d}\right)^2} {\displaystyle(-1)^{M}\prod_{\ell=1}^{k+M}\left(\ell\dfrac{\pi}{T_d}\right)\prod_{\ell=1}^{M-k}\left(-\ell\dfrac{\pi}{T_d}\right)}=\\
   &=& \lim_{M\rightarrow\infty}\dfrac{\displaystyle\prod_{m=1}^{M}\left(m\dfrac{\pi}{T_d}\right)^2}
   {\displaystyle(-1)^{M}(-1)^{M-k}\prod_{\ell=1}^{k+M}\left(\ell\dfrac{\pi}{T_d}\right)\prod_{\ell=1}^{M-k}\left(\ell\dfrac{\pi}{T_d}\right)}\\
   &=&\lim_{M\rightarrow\infty}\dfrac{\displaystyle\prod_{\ell=1}^{M}\left(\ell\dfrac{\pi}{T_d}\right)\prod_{\ell=1}^{M}\left(\ell\dfrac{\pi}{T_d}\right)} {\displaystyle (-1)^{2M-k}\prod_{\ell=1}^{k+M}\left(\ell\dfrac{\pi}{T_d}\right)\prod_{\ell=1}^{M-k}\left(\ell\dfrac{\pi}{T_d}\right)}=
   \lim_{M\rightarrow\infty}\dfrac{\displaystyle\prod_{\ell=M-k+1}^{M}\left(\ell\dfrac{\pi}{T_d}\right)} {\displaystyle (-1)^{-k}\prod_{\ell=M+1}^{M+k}\left(\ell\dfrac{\pi}{T_d}\right)}\\
   &=&\lim_{M\rightarrow\infty}(-1)^k\dfrac{\left[M-(k-1)\right]\left[M-(k-2)\right]\dots(M-1)M}{(M+1)(M+2)\dots\left[M+(k-1)\right](M+k)}\\
   &=&\lim_{M\rightarrow\infty}(-1)^k \dfrac{M-(k-1)}{M+k}\dfrac{M-(k-2)}{M+(k-1)}\dots\dfrac{M-1}{M+2}\dfrac{M}{M+1}\\
   &=&(-1)^k\lim_{M\rightarrow\infty}\dfrac{M-(k-1)}{M+k}\lim_{M\rightarrow\infty}\dfrac{M-(k-2)}{M+(k-1)}\dots
   \lim_{M\rightarrow\infty}\dfrac{M-1}{M+2}\lim_{M\rightarrow\infty}\dfrac{M}{M+1}=\\
   &=&(-1)^k\prod_{\ell=1}^{k}\lim_{M\rightarrow\infty}\dfrac{M-(\ell-1)}{M+\ell}
\end{eqnarray*}
Since we can always make $M$ infinitely greater than $k$ then
\begin{equation*}
  \lim_{M\rightarrow\infty}\dfrac{M-(\ell-1)}{M+\ell}=1,\quad\forall\ell=1,\dots,k
\end{equation*}
and, consequently,
\begin{equation*}
   b_k=(-1)^k=2T_da_k
\end{equation*}
and we conclude that
\begin{equation*}
   g(s)=\dfrac{\displaystyle\prod_{k=-\infty,k\neq0}^{\infty}\left(-jk\dfrac{\pi}{T_d}\right)}{\displaystyle\prod_{k=-\infty}^{\infty}(s-jk\dfrac{\pi}{T_d})}=
   2T_d\underbrace{\dfrac{e^{-T_ds}}{1-e^{-2T_ds}}}_{:=f(s)}
\end{equation*}
and the proof that the expansion of function $f(s)$ in Laurent series is possible is done in Appendix~\ref{RationalExpansionMF}.

$\hfill\Box$

\begin{theorem}\label{th3D}
The function
\begin{equation}\label{auxv}
v(s)=\dfrac{e^{-sT_d}}{1+e^{-2sTd}}
\end{equation}
can expanded as
\begin{equation*}
       v(s)=\sum_{k=-\infty}^{\infty}\dfrac{c_k}{s-\lambda_k}
\end{equation*}
where $c_k$ is the residual of $v(s)$ at $s=\lambda_k:$
$$ c_{k}=\dfrac{j(-1)^{(k+1)}}{2T_d}$$
\end{theorem}
\textbf{Proof:}
This function has poles at
\begin{equation*}
   \lambda_k=j(2k-1)\dfrac{\pi}{2T_d},\quad k=-\infty,\dots,-1,0,1,\infty
\end{equation*}
We can prove that condition \eqref{c1} holds for $v(s)$ exactly the same way we did for $f(s)$. So, we can expand $v(s)$ as
\begin{equation*}
       v(s)=\sum_{k=-\infty}^{\infty}\dfrac{c_k}{s-\lambda_k}
\end{equation*}
where $c_k$ is the residual of $v(s)$ at $s=\lambda_k$, i. e.
\begin{eqnarray*}
   c_k&=&\dfrac{1}{2j\pi}\oint_{\mathcal{C}_k}v(s)ds=\lim_{s\rightarrow \lambda_k}(s-\lambda_k)v(s)=
   \lim_{s\rightarrow \lambda_k}\dfrac{(s-\lambda_k)e^{-T_ds}}{1+e^{-2T_ds}}\\
   &=&
   \lim_{s\rightarrow \lambda_k}\dfrac{e^{-sT_d}-(s-\lambda_k)T_de^{-sT_d}}{-2T_de^{-2sT_d}}=\dfrac{e^{\lambda_kT_d}}{2T_d}=\dfrac{e^{\frac{j(2k-1)\pi}{2T_d}T_d}}{2Td}=
   \dfrac{e^{j(2k+1)\frac{\pi}{2}}}{2T_d}=\dfrac{j(-1)^{(k+1)}}{2T_d}
\end{eqnarray*}

$\hfill\Box$

\begin{theorem}\label{th4D}
The function
\begin{equation}\label{auxw}
w(s)=\dfrac{\displaystyle\prod_{k=-\infty,k\neq0}^{\infty}\left(-j(2k-1)\dfrac{\pi}{2T_d}\right)}
{\displaystyle\prod_{k=-\infty}^{\infty}\left(s-j(2k-1)\dfrac{\pi}{2T_d}\right)}
\end{equation}
can be written as:
$$w(s)=2\dfrac{e^{-T_ds}}{1+e^{-2T_ds}}=2v(s).$$
\end{theorem}
\textbf{Proof:}

Since $w(s)$ is  proper it can be expanded as
\begin{equation*}
   w(s)=\sum_{k=-\infty}^{\infty}\dfrac{d_k}{s-\lambda_k}
\end{equation*}
where $d_k$ is the residual of $w(s)$ at $s=\lambda_k=j(2k-1)\dfrac{\pi}{2T_d}$. In order to compute this residual we rewrite $w(s)$ as%\footnote{Notice that in order to keep real coefficients in the denominator polynomial complex poles must appear in conjugate pairs. This wouldn't happen if the inferior limit of the product was $-M$. For instance the  conjugate of $j(2k-1)\dfrac{\pi}{2T_d}$ is $j(-2k-3)\dfrac{\pi}{2T_d}$, not $j(-2k-1)\dfrac{\pi}{2T_d}$ (for $k=1 $ is the pole generated with $k=0$, not with $k=-1$).}}
\begin{equation*}
  w(s)=\lim_{M\rightarrow\infty}
  \dfrac{\displaystyle\prod_{k=-M+1}^{M}\left(-j(2k-1)\dfrac{\pi}{2T_d}\right)}{\displaystyle\prod_{k=-M+1}^{M}\left(s-j(2k-1)\dfrac{\pi}{2T_d}\right)}.
\end{equation*}
The residual $d_k$ is then given by
\begin{equation*}
   d_k=\dfrac{1}{2j\pi}\oint_{\mathcal{C}_k}w(s)ds=\lim_{s\rightarrow \lambda_k }(s-\lambda_k)w(s)= \lim_{M\rightarrow\infty}\dfrac{\displaystyle\prod_{m=-M+1}^{M}\left(j(2m-1)\dfrac{\pi}{2T_d}\right)} {\displaystyle\prod_{m=-M+1}^{k-1}\left(j(k-m)\dfrac{\pi}{T_d}\right)\prod_{m=k+1}^{M}\left(j(k-m)\dfrac{\pi}{T_d}\right)}
\end{equation*}
Given that the numerator can be expressed as the product of two complex factors we can write $d_k$ as
\begin{equation*}
   d_k=\lim_{M\rightarrow\infty}\dfrac{\displaystyle\prod_{m=1}^{M}\left((2m-1)\dfrac{\pi}{2T_d}\right)^2} {\displaystyle\prod_{m=-M+1}^{k-1}\left(j(k-m)\dfrac{\pi}{T_d}\right)\prod_{m=k+1}^{M}\left(j(k-m)\dfrac{\pi}{T_d}\right)}.
\end{equation*}
If now we replace define $\ell=k-m$ we get
\begin{eqnarray*}
   d_k &= & \lim_{M\rightarrow\infty}\dfrac{\displaystyle\prod_{m=1}^{M}\left((2m-1)\dfrac{\pi}{2T_d}\right)^2} {\displaystyle\prod_{\ell=1}^{k+M-1}\left(j\ell\dfrac{\pi}{T_d}\right)\prod_{\ell=-(M-k)}^{-1}\left(j\ell\dfrac{\pi}{T_d}\right)}=
   \lim_{M\rightarrow\infty}\dfrac{\displaystyle\prod_{m=1}^{M}\left((2m-1)\dfrac{\pi}{2T_d}\right)^2} {\displaystyle j^{2M-1}\prod_{\ell=1}^{k+M-1}\left(\ell\dfrac{\pi}{T_d}\right)\prod_{\ell=-(M-k)}^{-1}\left(\ell\dfrac{\pi}{T_d}\right)}=\\
   &=& \lim_{M\rightarrow\infty}\dfrac{\displaystyle\prod_{m=1}^{M}\left((2m-1)\dfrac{\pi}{2T_d}\right)^2} {\displaystyle(-j)^{2M}j^{-1}\prod_{\ell=1}^{k+M-1}\left(\ell\dfrac{\pi}{T_d}\right)\prod_{\ell=1}^{M-k}\left(-\ell\dfrac{\pi}{T_d}\right)}=\\
   &=&{ \lim_{M\rightarrow\infty}\dfrac{\displaystyle j\prod_{m=1}^{M}\left((2m-1)\dfrac{\pi}{2T_d}\right)^2}
   {\displaystyle(-1)^{M}(-1)^{M-k}\prod_{\ell=1}^{k+M-1}\left(\ell\dfrac{\pi}{T_d}\right)\prod_{\ell=1}^{M-k}\left(\ell\dfrac{\pi}{T_d}\right)}}=
   \dfrac{-j}{(-1)^{2M-k}T_d}=\dfrac{j(-1)^{k+1}}{T_d}
\end{eqnarray*}
and we conclude that
\begin{equation*}
   w(s)=\dfrac{\displaystyle\prod_{k=-\infty,k\neq0}^{\infty}\left(-j(2k+1)\dfrac{\pi}{2T_d}\right)} {\displaystyle\prod_{k=-\infty}^{\infty}(s-j(2k+1)\dfrac{\pi}{T_d})}=
   2\dfrac{e^{-T_ds}}{1+e^{-2T_ds}}=2v(s)
\end{equation*}

$\hfill\Box$

\begin{theorem}\label{Th5}
Consider the following functions
\begin{eqnarray}
\label{Th5-1}
G_e(\mathrm{S}) &=& \dfrac{\displaystyle\prod_{k=-\infty,k\neq0}^{\infty}\left(-jk\dfrac{\pi}{T_d}\right)}{\displaystyle\prod_{k=-\infty}^{\infty}(\mathrm{S}-jk\dfrac{\pi}{T_d})} \\
\label{Th5-2}
G_o(\mathrm{S})&=& \dfrac{\displaystyle\prod_{k=-\infty}^{\infty}\left(-j(2k-1)\dfrac{\pi}{2T_d}\right)}
{\displaystyle\prod_{k=-\infty}^{\infty}\left(\mathrm{S}-j(2k-1)\dfrac{\pi}{2T_d}\right)} \\
\label{Th5-3}
F_e(s) &=& \dfrac{e^{-T_d\mathrm{S}}}{1-e^{-2T_d\mathrm{S}}} \\
\label{Th5-4}
F_o(s) &=& \dfrac{e^{-\mathrm{S}T_d}}{1+e^{-2\mathrm{S}Td}}
\end{eqnarray}
then
\begin{eqnarray}
\label{Th5-5}
G_e(\mathrm{S}) &=& 2T_dF_e(\mathrm{S}) \\
\label{Th5-6}
G_o(\mathrm{S}) &=& 2F_o(\mathrm{S})
\end{eqnarray}
\end{theorem}
\textbf{Proof:}
To prove the results in Theorem~\ref{th2D} we consider $g(s)=G_{e}(s)$ and $f(s)=F_{e}(S)$ and also in  Theorem~\ref{th4D} we consider $w(s)=G_{o}(s)$ and $v(s)=F_{o}(S).$

$\hfill\Box$

$G_e(\mathrm{S})$ and $G_o(\mathrm{S})$ may also be written as\\\begin{eqnarray}
G_e(\mathrm{S}) &=& \dfrac{\displaystyle \prod_{k=1}^{\infty} 4k^2 \omega_0^2}{\displaystyle\prod_{k=0}^{\infty}\left(\mathrm{S}^2+4k^2\omega_0^2\right)}\\
G_o(\mathrm{S}) &=& \dfrac{ \displaystyle\prod_{k=1}^{\infty} (2k-1)^2 \omega_{0}^{2}}{
\displaystyle\prod_{k=1}^{\infty} \left(\mathrm{S}^2+(2k-1)^2\omega_0^2\right) }
\end{eqnarray}
where $\omega_o$ is defined in equation \eqref{30}.
As a result,
\begin{eqnarray*}
\displaystyle\prod_{k=1}^{\infty}\dfrac{\mathrm{S}^2+(2k-1)^2\omega_0^2}{\mathrm{S}^2+4k^2\omega_0^2}&=& \dfrac{\mathrm{S}G_e(\mathrm{S})}{G_o(\mathrm{S})} \displaystyle \prod_{k=1}^{\infty}\dfrac{1}{K_k} \\
\displaystyle\prod_{k=1}^{\infty}\dfrac{1}{\mathrm{S}^2+4k^2\omega_0^2}&=& \mathrm{S}G_e(\mathrm{S})\prod_{K=1}^{\infty}\dfrac{1}{4k^2\omega_0^2}
\end{eqnarray*}
with $K_k$ being defined in equation \eqref{31}. Now, using these equations and theorem \ref{Th5} we can write
\begin{equation}\label{35}
  \begin{array}{rcl}
     \prod_{k=1}^{\infty}\widehat{K}_k\dfrac{\mathrm{S}^2+(2k-1)^2\omega_0^2}{\mathrm{S}^2+4k^2\omega_0^2}&=& \bar{K}_{11}\dfrac{\mathrm{S}G_e(\mathrm{S})}{G_o(\mathrm{S})}=\\
     &=&\bar{K}_{11}T_d\dfrac{\mathrm{S}F_e(\mathrm{S})}{F_o(\mathrm{S})}\\
     \prod_{k=1}^{\infty} \dfrac{\alpha^2+4k^2\omega_0^2}{\mathrm{S}^2+4k^2\omega_0^2}=\bar{K}_{12}F_e(\mathrm{S})&=&2\bar{K}_{12}T_dG_e(\mathrm{S})
  \end{array}
\end{equation}
with
\begin{equation}\label{36}
  \begin{array}{l}
     \bar{K}_{11}=\displaystyle\prod_{k=1}^{\infty}\dfrac{\widehat{K}_k}{K_k}\\[5mm]
     \bar{K}_{12}=\displaystyle\prod_{k=1}^{\infty}\dfrac{\alpha^2+4k^2\omega_0^2}{4k^2\omega_0^2}
  \end{array}
\end{equation}
we thus can rewrite $\bar{G}_{11}(\mathrm{S})$ and $\bar{G}_{21}(\mathrm{S})$ as
\begin{equation}\label{37}
  \begin{array}{l}
      \bar{G}_{11}(\mathrm{S})=\dfrac{K_GT_d\bar{K}_{11}\mathrm{S}\left(1+e^{-2\mathrm{S}T_d}\right)} {\left(\mathrm{S}-\alpha\right)\left(1-e^{-2\mathrm{S}T_d}\right)}\\[5mm]
      \bar{G}_{21}(\mathrm{S})=\dfrac{2K_GT_d\bar{K}_{12}\mathrm{S}e^{-\mathrm{S}T_d}}{\left(\mathrm{S}-\alpha\right)\left(1-e^{-2\mathrm{S}T_d}\right)}.
  \end{array}
\end{equation}
We can compute $\widehat{G}_{11}(s)$ and $ \widehat{G}_{21}(s)$ from these equations by replacing $\mathrm{S}$ with $s+\alpha$:
\begin{equation}\label{38}
  \begin{array}{l}
      \widehat{G}_{11}(s)=\dfrac{K_GT_d\bar{K}_{11}\left(s+\alpha\right)\left(1+e^{-2\alpha T_d}e^{-2sT_d}\right)} {s\left(1-e^{-2\alpha T_d} e^{-2sT_d}\right)}\\[5mm]
      \widehat{G}_{21}(s)=\dfrac{2K_GT_d\bar{K}_{12}\left(s+\alpha\right)e^{-\alpha T_d} e^{-sT_d}}{s\left(1-e^{-2\alpha T_d}e^{-2sT_d}\right)}.
  \end{array}
\end{equation}
Given that $\widehat{G}_{11}(s)$ and $\widehat{G}_{21}(s)$ were defined in such a way that
\begin{equation}\label{39}
  \begin{array}{rcl}
     \displaystyle
     \lim_{s\rightarrow 0}s\widehat{G}_{11}(s) & = &
     \displaystyle\lim_{s\rightarrow 0}s\widehat{G}_{21}(s)=\lim_{s\rightarrow 0}sG_{11}(s)\\
     & = & \displaystyle\lim_{s\rightarrow 0}sG_{21}(s)=K_G,
  \end{array}
\end{equation}
we can rewrite \eqref{38} as
\begin{equation}\label{40}
  \begin{array}{l}
      \widehat{G}_{11}(s)=K_{11}\dfrac{\left(s+\alpha\right)\left(1+e^{-2\alpha T_d}e^{-2sT_d}\right)} {s\left(1-e^{-2\alpha T_d} e^{-2sT_d}\right)}\\[5mm]
      \widehat{G}_{21}(s)=K_{21}\dfrac{\left(s+\alpha\right) e^{-sT_d}}{s\left(1-e^{-2\alpha T_d}e^{-2sT_d}\right)}.
  \end{array}
\end{equation}
with
\begin{equation}\label{41}
  \begin{array}{l}
     K_{11}=\dfrac{K_G\left(1-e^{-2\alpha T_d}\right)}{\alpha \left(1+e^{-2\alpha T_d}\right)}\\[5mm]
     K_{21}=\dfrac{K_G\left(1-e^{-2\alpha T_d}\right)}{\alpha}.
  \end{array}
\end{equation}
\section{Case study}\label{CaseStudy}
In this section we study a small pipeline using the lumped linear model derived above. The pipe has a length $L=35\text{Km}$ and a diameter $\mathcal{D}=793\text{mm}$. The friction factor is $f_c=0.0079$ and the isothermal speed of sound is $c=300\text{m/sec}$. We considered  $q_m=90{Kg/sec}$ and  $p_m=80\text{ bar}=8\times 10^6\text{ Pascal}$ as the mass-flow and pressure nominal values. The $\alpha$ and $K_G$ parameters were calculated from equations \eqref{2} and Corollary~\ref{cor5}, respectively, and the values of $0.0051$ and $5.064$ were obtained. We computed the frequency responses (FR) of $G_{11}(s)$ and $G_{21}(s)$   from equations \eqref{29} using truncated approximations of order $n$, $G_{11}^{(n)}(s)$ and $G_{11}^{(n)}(s)$, respectively, in a frequency bandwidth, $BW=\left[10^{-4}\text{rad/sec},\;1\;\text{rad/sec}\right]$.
\begin{figure}[h]
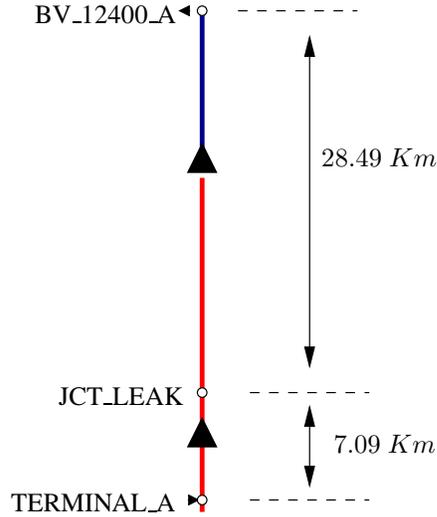

\begin{center}
\scalebox{1}{ \hspace{0cm} \input tipologia.pstex_t }\\
\caption{Gas pipeline topology.}
\label{fig1}
\end{center}
\end{figure}

\begin{figure}
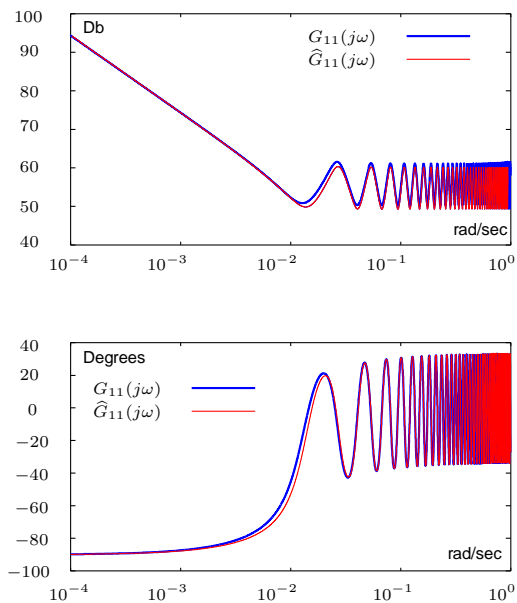

\begin{center}
\scalebox{0.90}{ \hspace{0cm} \input G11.pstex_t }
\caption{Bode plots of  $G_{11}(s)$ and the approximation $\widehat{G}_{11}(s)$.}
\label{f1}
\end{center}
\end{figure}

\begin{figure}
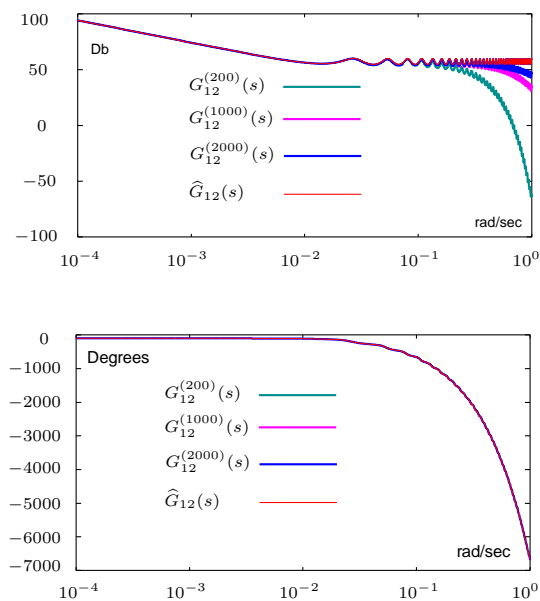

\begin{center}
\scalebox{0.90}{ \hspace{0cm} \input G12.pstex_t }
\caption{Bode plots of both $G_{12}^{(n)}(s)$, $n=200,\;100,\;2000$ and the approximation $\widehat{G}_{12}(s)$.}
\label{fig2}
\end{center}
\end{figure}

Figures \ref{f1} and \ref{fig2} display the respective Bode diagrams and compare them with their approximations $\widehat{G}_{11}(s)$ and $\widehat{G}_{21}(s)$. $G_{11}(s)$ converges very fast. The FR of a truncated approximation with two hundred poles and zeros  had already converged to its limit in the  whole frequency interval $BW$. Figure \ref{f1}  shows that there are no significant differences between $G_{11}(s)$ and $\widehat{G}_{11}(s)$. Consequently, $G_{11}(s)$ can be substituted by its approximation without loss of accuracy and with a significant reduction of the computational costs.

The convergence of $G_{21}(s)$ is much slower. A two thousand order approximation didn't converge in the whole bandwidth $BW$. We can subsequently conclude that a truncated approximation of equation \eqref{29} needs too many factors leading, therefore,  to high order transfer functions with high computational costs.

\begin{figure}[h]
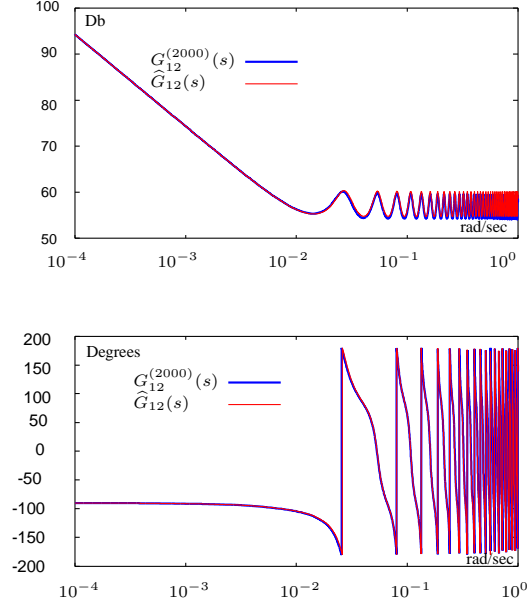

\begin{center}
\scalebox{0.90}{ \hspace{0cm} \input G12-Novo.pstex_t }
\caption{Higher resolution Bode plots of both $G_{12}^{(2000)}(s)$ and the approximation $\widehat{G}_{12}(s)$.}
\label{f3}
\end{center}
\end{figure}

\begin{figure}
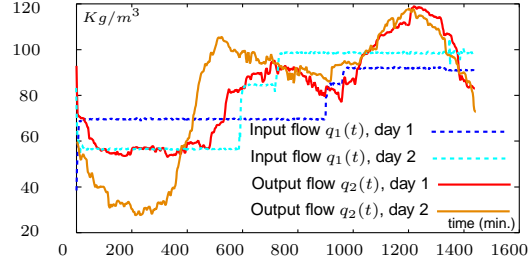

\begin{center}
\scalebox{0.90}{ \hspace{0cm} \input QiQo.pstex_t }
\caption{Input and output mass-flows on the both days.}
\label{f4}
\end{center}
\end{figure}

However,   from Figures~\ref{fig2} and \ref{f3} one can conclude that  $G_{21}(s)$ can be also substituted by its approximation. In Figure \ref{f3} the Bode diagrams of $\widehat{G}_{21}(s)$ and $G_{21}^{(2000)}(s)$ are compared with a better resolution, i.e.  the phase is restricted to the range $(-180 \text{ degrees},\;180 \text{ degrees}]$. The phase discontinuities are due to the phase-crossing of the odd multiples of 180  degrees which are converted from -180 degrees to 180 degrees. Also notice that the Bode diagram of the finite approximation $G_{12}^{(2000)}(s)$ converges to $G_{12}(s)$ almost in the whole frequency $BW$.

\begin{figure}
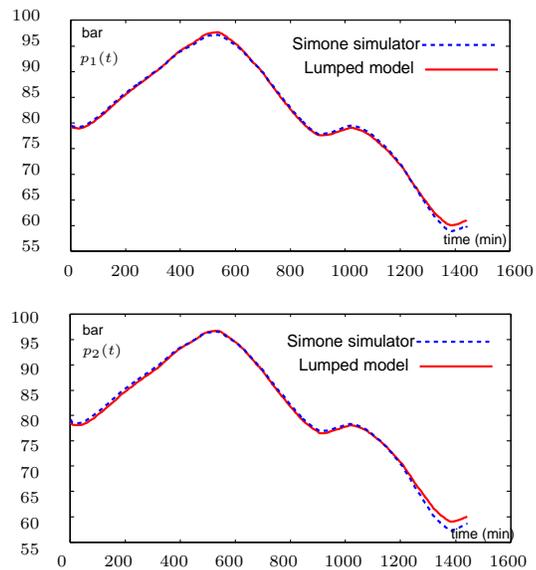

\begin{center}
\scalebox{0.90}{ \hspace{0cm} \input Psim-Dia1.pstex_t }
\caption{Input and Output pressures, $P_1$, and $P_2$, respectively, simulated by the lumped model \eqref{26} (solide red line) and by SIMONE$^{\textregistered }$ using the first day data.}
\label{f6}
\end{center}
\end{figure}

This pipeline was simulated taking two normal days operation data as the input and output mass-flows. The simulation was performed with the previous referred
SIMONE$^{\textregistered }$  simulator. Figure \ref{f4} shows the input and output mass-flows on both days

The intake and offtake massflows are denoted by $q_1(t)$ and $q_2(t)$, respectively.

\begin{figure}
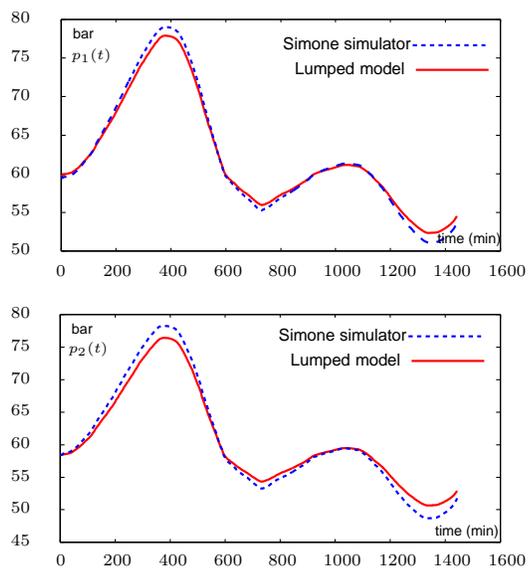

\begin{center}
\scalebox{0.90}{ \hspace{0cm} \input Psim-Dia2.pstex_t }
\caption{Input and Output pressures, $P_1$, and $P_2$, respectively, simulated by the lumped model \eqref{26} (solide red line) and by SIMONE$^{\textregistered }$ using the second day data.}
\label{f7}
\end{center}
\end{figure}
Figures \ref{f6} and \ref{f7} compare the input and output pressures simulated with the lumped transfer function model
\begin{equation}\label{26}
 \begin{array}{rcl}
    P_{1}(s) &=& \widehat{G}_{11}(s)Q_1(s) - \widehat{G}_{21}(s)Q_2(s) \\
    P_{2}(s) &=& \widehat{G}_{21}(s)Q_1(s) - \widehat{G}_{22}(s)Q_2(s)
 \end{array}
\end{equation}
with the ones simulated with SIMONE$^{\textregistered }$  for the two days. The intake and offtake pressures are denoted by $p_1(t)$ and $p_2(t)$, respectively. We can see that the results are better for the first day. This is expectable, since the mass-flows and pressures of the second day data present stronger deviations from the  nominal values used for the calculation of the $\alpha$ parameter.
As a matter of fact, on the second day the quotient $\left(q_1(t)+q_2(t)\right)/\left(p_1(t)+p_2(t)\right)$  has  a root mean square deviation from the nominal value of about $50\%$ (this deviation is of $30\%$ in the first day). Yet, in both cases, the model has well captured the dynamics of the system and, therefore, it seems to be a valuable tool for gas leakage detection and gas networks controller design.

%==========================================
\appendix
%===========================================
\section{Change of variables in an integral model for a short gas pipeline}\label{changeofvariables}
In this section, a change in the state-space variables of the integral model is performed. The purpose is to obtain a simpler system matrix in order to simplify the determination of the eigenvalues of the system.

Consider model \eqref{10}  with the respective matrices defined according to \eqref{12}--\eqref{17}.

Consider for this system the following change of variables:
\bean
\label{A-1}
\ba{rcl}
z_2&=&x_2+x_2+\cdots+x_N+x_{N+1}\\
z_2&=&x_2-x_2\\
z_3&=&x_2-x_3\\
&\vdots&\\
z_{N+1}&=&x_{N}-x_{N+1}\\
z_{N+2}&=&x_{N+2}\\
z_{N+3}&=&x_{N+3}\\
& \vdots &\\
z_{2N+1}&=&x_{2N+1}
\ea
\eean
and we obtain the following realisation
\bean
\label{A-2}
\ba{rcl}
\dot{z}_2(t)&=&\dfrac{c^2}{\mathcal{A} \Delta \ell}u_2(t)-\dfrac{c^2}{\mathcal{A} \Delta \ell}u_2(t)\\[4mm]
\dot{z}_2(t)&=&-2\dfrac{c^2}{\mathcal{A} \Delta \ell}z_{N+2}(t)+\dfrac{c^2}{\mathcal{A} \Delta \ell}z_{N+3}(t)\\[4mm]
\dot{z}_3(t)&=&\dfrac{c^2}{\mathcal{A} \Delta \ell}z_{N+2}(t)-2\dfrac{c^2}{\mathcal{A} \Delta \ell}z_{N+3}(t)
+\dfrac{c^2}{\mathcal{A} \Delta \ell}z_{N+4}(t)\\[4mm]
&\vdots&\\[4mm]
\dot{z}_i(t)&=&\dfrac{c^2}{\mathcal{A} \Delta \ell}z_{N+i-1}(t)-2\dfrac{c^2}{\mathcal{A} \Delta \ell}z_{N+i}(t)
+\dfrac{c^2}{\mathcal{A} \Delta \ell}z_{N+i+1}(t)
\ea
\eean
\bean
\label{A-2a}
\ba{rcl}
%\\[4mm]
\dot{z}_{N+1}(t)&=&\dfrac{c^2}{\mathcal{A} \Delta \ell}z_{2N}(t)-2\dfrac{c^2}{\mathcal{A} \Delta \ell}z_{2N+1}(t)
+\dfrac{c^2}{\mathcal{A} \Delta \ell}u_2(t)\\[4mm]
\dot{z}_{N+2}(t)&=&\dfrac{\mathcal{A}}{\Delta \ell}z_2(t)-
\dfrac{f_c c^2 Q_{m1}}{2 \mathcal{D} \mathcal{A} P_{m1}}z_{N+2}(t)\\[4mm]
\dot{z}_{N+3}(t)&=&\dfrac{\mathcal{A}}{\Delta \ell}z_3(t)-
\dfrac{f_c c^2 Q_{m1}}{2 \mathcal{D} \mathcal{A} P_{m1}}z_{N+3}(t)\\[4mm]
&\vdots&
\ea
\eean
\bean
\label{A-2a}
\ba{rcl}
%\\[4mm]
\dot{z}_{N+i}(t)&=&\dfrac{\mathcal{A}}{\Delta \ell}z_i(t)-
\dfrac{f_c c^2 Q_{m1}}{2 \mathcal{D} \mathcal{A} P_{m1}}z_{N+i}(t)\\[4mm]
&\vdots&\\[4mm]
\dot{z}_{2N+1}(t)&=&\dfrac{\mathcal{A}}{\Delta \ell}z_{N+1}(t)-
\dfrac{f_c c^2 Q_{m1}}{2 \mathcal{D} \mathcal{A} P_{m1}}z_{2N+1}(t)
\ea
\eean
\bean
\label{A-2a}
\ba{rcl}
y_2(t)&=&\dfrac{1}{N+1}z_2(t)+\dfrac{N}{N+1}z_2(t)+\dfrac{N-1}{N+1}z_3(t)+\dots+\dfrac{1}{N+1}z_{N+1}(t)\\[4mm]
y_2(t)&=&\dfrac{1}{N+1}z_2(t)-\dfrac{1}{N+1}z_2(t)-\dfrac{2}{N+1}z_3(t)-\dots-\dfrac{N}{N+1}z_{N+1}(t)
\ea
\eean
and in matricial form, we have:
\bean
\label{A-3}
\ba{rcl}
\dot{z}(t)&=&
\bar{A}z(t)+\bar{B}u(t)\\
y(t)&=&\bar{C}z(t)
\ea
\eean
 with
\bean
\label{A-4}
\ba{rcl}
\bar{A}&=&\left[
\ba{c|c|c}
0 & 0_{1\times N} & 0_{1\times N} \\
\hline
& & \\
0_{N \times 1} &  \bar{A}_{11} & \bar{A}_{12} \\
\hline
& & \\
0_{N \times 1} & \bar{A}_{21} & \bar{A}_{22}\\
\ea
\right]
\ea
\eean
\bean
\label{A-2b}
\ba{rcl}
\bar{A}_{11}&=& 0_{N \times N}
\ea
\eean
\bean
\label{A-2b}
\ba{rcl}
\bar{A}_{12}&=&\left[
\ba{ccccccc}
-2\dfrac{c^2}{\mathcal{A} \Delta \ell} & \dfrac{c^2}{\mathcal{A} \Delta \ell} & 0 & \cdots & 0 & 0& 0\\
[4mm]
\dfrac{c^2}{\mathcal{A} \Delta \ell} & -2\dfrac{c^2}{\mathcal{A} \Delta \ell} & \dfrac{c^2}{\mathcal{A} \Delta \ell} &
\cdots &  \vdots &  \vdots &  \vdots\\
[4mm]
0 &\dfrac{c^2}{\mathcal{A} \Delta \ell} & -2\dfrac{c^2}{\mathcal{A} \Delta \ell} & \ddots  &  \vdots  &  \vdots&  \vdots\\
\vdots & \vdots & \ddots &  \ddots & \ddots & \cdots & \cdots\\
 \vdots &  \vdots &  \vdots & \cdots & \dfrac{c^2}{\mathcal{A} \Delta \ell} & -2\dfrac{c^2}{\mathcal{A} \Delta \ell} &
\dfrac{c^2}{\mathcal{A} \Delta \ell}\\
[4mm]
0 & 0 & 0 & \cdots & 0 & \dfrac{c^2}{\mathcal{A} \Delta \ell} & -2\dfrac{c^2}{\mathcal{A} \Delta \ell}
\ea
\right]
\in \cR^{N \times N}
\ea
\eean
\bean
\label{A-2a}
\ba{rcl}
\bar{A}_{21}&=&\dfrac{\mathcal{A}}{\Delta \ell}I_N, \quad I_N\;-\text{identity matrix} \;N \times N\\[4mm]
\bar{A}_{22}&=&-\dfrac{f_c c^2 Q_{m}}{2 \mathcal{D} \mathcal{A} P_{m}}I_N\\[4mm]
\bar{B}&=& \dfrac{c^2}{\mathcal{A} \Delta \ell} \left[\ba{c|c} e_2 & -e_2+e_{N+1} \ea \right]\in \cR^{(2N+1) \times 2} \\[4mm]
\bar{C}&=&\left[\ba{ccccc} \dfrac{1}{N+1} & \dfrac{N}{N+1} & \dfrac{N-1}{N+1} & \cdots & \dfrac{1}{N+1} \\[4mm]
                     \dfrac{1}{N+1} & -\dfrac{1}{N+1} & -\dfrac{2}{N+1} & \cdots & -\dfrac{N}{N+1}
\ea \right]
\in \cR^{2 \times (2N+1)}
\ea
\eean
Matrices $A$ and $\bar{A}$ have the same spectrum, however its calculation seems to be easier if we use matrix $\bar{A}.$
%%======================================
\section{Rational expansion of meromorphic functions}\label{RationalExpansionMF}
Let $f(s)$ be a function meromorphic in the finite complex plane with poles at $\lambda_1,\;\lambda_2,\;\ldots$, and let ($\Gamma_1,\;\Gamma_2\;,\ldots$) be a sequence of simple closed curves such that:
\begin{itemize}
 \item The origin lies inside each curve $\Gamma_k$.
 \item No curve passes through a pole of $f$.
 \item $\Gamma_k$ lies inside $\Gamma_{k+1}$ for all k
 \item $\displaystyle\lim_{k\rightarrow \infty} d(\Gamma_k) = \infty$, where $d(\Gamma_k)$ gives the distance from the curve to the origin
\end{itemize}
Suppose also that there exists an integer p such that
\begin{equation*}
   \displaystyle
   \lim_{k\rightarrow \infty} \oint_{\Gamma_k} \left|\frac{f(s)}{s^{p}}\right|\left| ds\right| < \infty.
\end{equation*}
Denoting the principal part of the Laurent series of $f$ about the point $\lambda_k$ as $PP\left[f(s);s=\lambda_k\right]$, we have, if $p<0$.
\begin{equation*}
       f(z) = \sum_{k=0}^{\infty} \operatorname{PP}(f(z); z = \lambda_k).
\end{equation*}

\begin{theorem}\label{th1D}
Consider function $f(s)=\dfrac{e^{-T_ds}}{1-e^{-2T_ds}}.$  There exists an integer p such that
\begin{equation}
   \displaystyle
   \lim_{k\rightarrow \infty} \oint_{\Gamma_k} \left|\frac{f(s)}{s^{p}}\right|\left| ds\right| < \infty.
\end{equation}
\end{theorem}
\textbf{Proof:}
This function has poles at
\begin{equation*}
   \lambda_k=j2k\dfrac{\pi}{2T_d}=jk\dfrac{\pi}{T_d},\quad k=-\infty,\dots,-1,0,1,\infty
\end{equation*}
The contours $\Gamma_k$ will be squares vertices at $\pm(2k-1)\dfrac{\pi}{2T_d}\pm j(2k-1)\dfrac{\pi}{2T_d}$, $k>1$, traversed counterclockwise, which are easily seen to satisfy the necessary conditions.
\begin{figure}[ht]
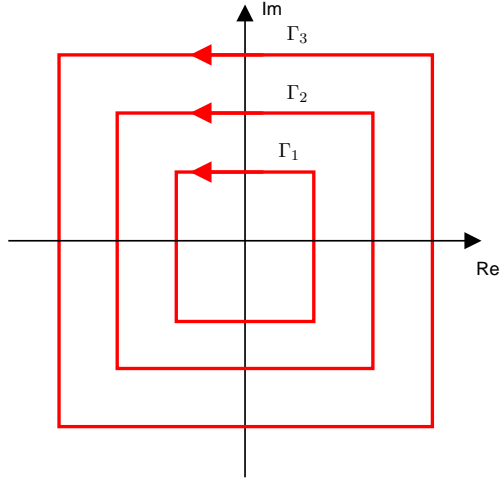

\begin{center}
\scalebox{0.7}{ \hspace{0cm} \input ContourSquares.pstex_t }
\caption{Contours}
\label{f2}
\end{center}
\end{figure}
To see what are the terms of the Laurent series expansion of $f(s)$ we need do see for which $p$ the condition
\begin{equation}\label{c1}
   \displaystyle
   \lim_{k\rightarrow \infty} \oint_{\Gamma_k} \left|\frac{f(s)}{s^{p}}\right| \left|ds\right| < \infty
\end{equation}
holds. We can partition this integral as
\begin{equation*}
   \displaystyle
   \oint_{\Gamma_k} \left|\frac{f(s)}{s^{p}}\right| \left|ds\right|=\oint_{\Gamma_{k1}}\left|\frac{f(s)}{s^{p}}\right| \left|ds\right|+
   \oint_{\Gamma_{k2}}\left|\frac{f(s)}{s^{p}}\right| \left|ds\right|+\oint_{\Gamma_{k3}}\left|\frac{f(s)}{s^{p}}\right| \left|ds\right|+
   \oint_{\Gamma_{k4}}\left|\frac{f(s)}{s^{p}}\right| \left|ds\right|
\end{equation*}
where
\begin{eqnarray*}
 \Gamma_{k1} &=& \left\{s:s=(2k-1)\dfrac{\pi}{2T_d}+j\omega,\; \omega \uparrow,\;\omega \in [-(2k-1)\dfrac{\pi}{2T_d},(2k-1)]\dfrac{\pi}{2T_d}\right\} \\
 \Gamma_{k2} &=& \left\{s:s=\sigma+j2(k-1)\dfrac{\pi}{2T_d},\; \sigma \downarrow,\;\sigma \in [-(2k-1)\dfrac{\pi}{2T_d},(2k-1)]\dfrac{\pi}{2T_d}\right\} \\
 \Gamma_{k3} &=& \left\{s:s=-(2k-1)\dfrac{\pi}{2T_d}+j\omega,\; \omega \downarrow,\;\omega \in [-(2k-1)\dfrac{\pi}{2T_d},(2k-1)]\dfrac{\pi}{2T_d}\right\} \\
 \Gamma_{k4} &=& \left\{s:s=\sigma-j(2k-1)\dfrac{\pi}{2T_d},\; \sigma \uparrow,\;\omega \in [-(2k-1)\dfrac{\pi}{2T_d},(2k-1)]\dfrac{\pi}{2T_d}\right\}.
\end{eqnarray*}
Notice that
\begin{itemize}
  \item For $s\in\Gamma_{k1}$, $|ds|=d\omega$.
  \item For $s\in\Gamma_{k2}$, $|ds|=-d\sigma$.
  \item For $s\in\Gamma_{k3}$, $|ds|=-d\omega$.
  \item For $s\in\Gamma_{k4}$, $|ds|=d\sigma$.
\end{itemize}
Now we can write
\begin{equation*}
  \begin{array}{l}
   \displaystyle
   \lim_{k\rightarrow \infty} \oint_{\Gamma_k} \left|\frac{f(s)}{s^{p}}\right| \left|ds\right| =\displaystyle\int_{-(2k-1)\frac{\pi}{2T_d}}^{(2k-1)\frac{\pi}{2T_d}}\left|f\left((2k-1)\frac{\pi}{2T_d}+j\omega\right)\right|d\omega -\displaystyle\int_{(2k-1)\frac{\pi}{2T_d}}^{-(2k-1)\frac{\pi}{2T_d}}\left|f\left(\sigma+j(2k-1)\frac{\pi}{2T_d}\right)\right|d\sigma \\ [5mm] -\displaystyle\int_{(2k-1)\frac{\pi}{2T_d}}^{-(2k-1)\frac{\pi}{2T_d}}\left|f\left(-(2k-1)\frac{\pi}{2T_d}+j\omega\right)\right|d\omega +\displaystyle\int_{-(2k-1)\frac{\pi}{2T_d}}^{(2k-1)\frac{\pi}{2T_d}}\left|f\left(\sigma-j(2k-1)\frac{\pi}{2T_d}\right)\right|d\sigma\\[5mm]
   =\displaystyle\int_{-(2k-1)\frac{\pi}{2T_d}}^{(2k-1)\frac{\pi}{2T_d}}\left|f\left((2k-1)\frac{\pi}{2T_d}+j\omega\right)\right|d\omega +\displaystyle\int_{-(2k-1)\frac{\pi}{2T_d}}^{(2k-1)\frac{\pi}{2T_d}}\left|f\left(\sigma+j(2k-1)\frac{\pi}{2T_d}\right)\right|d\sigma \\ [5mm] +\displaystyle\int_{-(2k-1)\frac{\pi}{2T_d}}^{(2k-1)\frac{\pi}{2T_d}}\left|f\left(-(2k-1)\frac{\pi}{2T_d}+j\omega\right)\right|d\omega +\displaystyle\int_{-(2k-1)\frac{\pi}{2T_d}}^{(2k-1)\frac{\pi}{2T_d}}\left|f\left(\sigma-j(2k-1)\frac{\pi}{2T_d}\right)\right|d\sigma\\[5mm]
  \end{array}
\end{equation*}
Next we analyze the four terms of this integral
\begin{description}
\item[First term]
\begin{equation*}
   \begin{array}{c}
       \displaystyle\int_{-(2k-1)\frac{\pi}{2T_d}}^{(2k-1)\frac{\pi}{2T_d}}\left|f\left((2k-1)\frac{\pi}{2T_d}+j\omega\right)\right|d\omega =
       \displaystyle\int_{-(2k-1)\frac{\pi}{2T_d}}^{(2k-1)\frac{\pi}{2T_d}}\left|\dfrac{e^{-(2k-1)\frac{\pi}{2}}e^{-j\omega T_d}} {\left((2k-1)\frac{\pi}{2T_d}+j\omega\right)^p\left(1-e^{-(2k-1)\pi}e^{-j\omega T_d}\right)}\right|d\omega=\\[5mm]
       \displaystyle\int_{-(2k-1)\frac{\pi}{2T_d}}^{(2k-1)\frac{\pi}{2T_d}}\left|\dfrac{e^{-(2k-1)\frac{\pi}{2}}} {\left((2k-1)\frac{\pi}{2T_d}+j\omega\right)^p\left(1-e^{-(2k-1)\pi}e^{-j\omega T_d}\right)}\right|d\omega.
   \end{array}
\end{equation*}
Given that
\begin{equation*}
   \lim_{k\rightarrow\infty}e^{-(2k-1)\frac{\pi}{2}}=0
\end{equation*}
then \small
\begin{equation*}
   \hspace{-1cm}\displaystyle\lim_{k\rightarrow\infty} \int_{-(2k-1)\frac{\pi}{2T_d}}^{(2k-1)\frac{\pi}{2T_d}}\left|f\left((2k-1)\frac{\pi}{2T_d}+j\omega\right)\right|d\omega=
   \displaystyle\lim_{k\rightarrow\infty}\int_{-(2k-1)\frac{\pi}{2T_d}}^{(2k-1)\frac{\pi}{2T_d}}\left|\dfrac{e^{-(2k-1)\frac{\pi}{2}}} {\left((2k-1)\frac{\pi}{2T_d}+j\omega\right)^p\left(1-e^{-(2k-1)\pi}e^{-j\omega T_d}\right)}\right|d\omega=0.
   \end{equation*}
   \normalsize
\item[Second term]
\small
\begin{equation*}
  \begin{array}{c}
   \displaystyle\int_{-(2k-1)\frac{\pi}{2T_d}}^{(2k-1)\frac{\pi}{2T_d}}\left|f\left(\sigma+j(2k-1)\frac{\pi}{2T_d}\right)\right|d\sigma
   =\displaystyle\int_{-(2k-1)\frac{\pi}{2T_d}}^{(2k-1)\frac{\pi}{2T_d}}\left|\dfrac{e^{-\sigma T_d}e^{-j(2k-1)\frac{\pi}{2}}} {\left(\sigma+j(2k-1)\frac{\pi}{2T_d}\right)^p\left(1-e^{-2\sigma T_d}e^{-j(2k-1)\pi}\right)}\right|d\sigma=\\[8mm]
   \hspace{-15mm}=\displaystyle\int_{-(2k-1)\frac{\pi}{2T_d}}^{0}\left|\dfrac{e^{-\sigma T_d}}
   {\left(\sigma+j(2k-1)\frac{\pi}{2T_d}\right)^p\left(1-e^{-2\sigma T_d}(-1)^k\right)}\right|d\sigma+
   \displaystyle\int_{0}^{(2k-1)\frac{\pi}{2T_d}}\left|\dfrac{e^{-\sigma T_d}}
   {\left(\sigma+j(2k-1)\frac{\pi}{2T_d}\right)^p\left(1-e^{-2\sigma T_d}(-1)^k\right)}\right|d\sigma=\\
  \end{array}
\end{equation*}
\normalsize
Notice that
\small
\begin{equation*}
   \displaystyle \int_{-(2k-1)\frac{\pi}{2T_d}}^{0}\left|\dfrac{e^{-\sigma T_d}}
   {\left(\sigma+j(2k-1)\frac{\pi}{2T_d}\right)^p\left(1-e^{-2\sigma T_d}(-1)^k\right)}\right|d\sigma \leq
    \displaystyle \int_{-(2k-1)\frac{\pi}{2T_d}}^{0}\left|\dfrac{e^{-\sigma T_d}}
   {\left(\sigma+j(2k-1)\frac{\pi}{2T_d}\right)^p\left(1-e^{-2\sigma T_d}\right)}\right|d\sigma
\end{equation*}
\normalsize
On the other hand, for $\sigma <0$
\begin{equation*}
   \left|1-e^{-\sigma T_d}\right|=e^{|\sigma|T_d}-1<e^{|\sigma T_d|}=e^{-\sigma T_d}.
\end{equation*}
As a result
\begin{equation*}
 \begin{array}{c}
   \displaystyle \int_{-(2k-1)\frac{\pi}{2T_d}}^{0}\left|\dfrac{e^{-\sigma T_d}}
   {\left(\sigma+j(2k-1)\frac{\pi}{2T_d}\right)^p\left(1-e^{-2\sigma T_d}\right)}\right|d\sigma <
   \displaystyle \int_{-(2k-1)\frac{\pi}{2T_d}}^{0}\left|\dfrac{e^{-\sigma T_d}}
   {\left(\sigma+j(2k-1)\frac{\pi}{2T_d}\right)^pe^{-2\sigma T_d}}\right|d\sigma=\\[8mm]
   =\displaystyle \int_{-(2k-1)\frac{\pi}{2T_d}}^{0}\left|\dfrac{e^{\sigma T_d}}
   {\left(\sigma+j(2k-1)\frac{\pi}{2T_d}\right)^p}\right|d\sigma
 \end{array}
\end{equation*}
Making $k\rightarrow\infty$,
\begin{equation*}
   \displaystyle\lim_{k\rightarrow \infty}\displaystyle \int_{-(2k-1)\frac{\pi}{2T_d}}^{0}\left|\dfrac{e^{\sigma T_d}}
   {\left(\sigma+j(2k-1)\frac{\pi}{2T_d}\right)^p}\right|d\sigma=M_1.\footnote{Remember that in this integral $\sigma$ only takes negative values.}
\end{equation*}
Notice also that
\small
\begin{equation*}
   \displaystyle\int_{0}^{(2k-1)\frac{\pi}{2T_d}}\left|\dfrac{e^{-\sigma T_d}}
   {\left(\sigma+j(2k-1)\frac{\pi}{2T_d}\right)^p\left(1-e^{-2\sigma T_d}(-1)^k\right)}\right|d\sigma\leq
   \displaystyle\int_{0}^{(2k-1)\frac{\pi}{2T_d}}\left|\dfrac{e^{-\sigma T_d}}
   {\left(\sigma+j(2k-1)\frac{\pi}{2T_d}\right)^p\left(1-e^{-2\sigma T_d}\right)}\right|d\sigma
\end{equation*}
\normalsize
For $\sigma >0$
\begin{equation*}
   1-e^{-\sigma Td}<1.
\end{equation*}
Consequently
\begin{equation*}
   \displaystyle\int_{0}^{(2k-1)\frac{\pi}{2T_d}}\left|\dfrac{e^{-\sigma T_d}}
   {\left(\sigma+j(2k-1)\frac{\pi}{2T_d}\right)^p\left(1-e^{-2\sigma T_d}\right)}\right|d\sigma <
   \displaystyle\int_{0}^{(2k-1)\frac{\pi}{2T_d}}\left|\dfrac{e^{-\sigma T_d}}
   {\left(\sigma+j(2k-1)\frac{\pi}{2T_d}\right)^p}\right|d\sigma
\end{equation*}
Making $k\rightarrow\infty$ again
\begin{equation*}
   \displaystyle\lim_{k\rightarrow\infty}\int_{0}^{(2k-1)\frac{\pi}{2T_d}}\left|\dfrac{e^{-\sigma T_d}}
   {\left(\sigma+j(2k-1)\frac{\pi}{2T_d}\right)^p}\right|d\sigma=M_1
\end{equation*}
and we conclude that
\begin{equation*}
   \displaystyle\lim_{k\rightarrow\infty} \int_{-(2k-1)\frac{\pi}{2T_d}}^{(2k-1)\frac{\pi}{2T_d}}\left|f\left(\sigma+j(2k-1)\frac{\pi}{2T_d}\right)\right|d\sigma<2M_1<\infty.
\end{equation*}
\item[Third term] This term is similar to the first one and using the same arguments we can prove that
\begin{equation*}
   \lim_{k\rightarrow\infty} \displaystyle\int_{-(2k-1)\frac{\pi}{2T_d}}^{(2k-1)\frac{\pi}{2T_d}}\left|f\left(-(2k-1)\frac{\pi}{2T_d}+j\omega\right)\right|d\omega=0.
\end{equation*}
\item In similar way that we did for the second term we can prove that
\begin{equation*}
   \displaystyle\int_{-(2k-1)\frac{\pi}{2T_d}}^{(2k-1)\frac{\pi}{2T_d}}\left|f\left(\sigma-j(2k-1)\frac{\pi}{2T_d}\right)\right|d\sigma<2M_1
\end{equation*}
\end{description}
We can now conclude that condition \eqref{c1} holds for any $p<0$.

$\hfill\Box$

Since there exists an integer p such that
\begin{equation*}
   \displaystyle
   \lim_{k\rightarrow \infty} \oint_{\Gamma_k} \left|\frac{f(s)}{s^{p}}\right|\left| ds\right| < \infty.
\end{equation*}
Then, denoting the principal part of the Laurent series of $f$ about the point $\lambda_k$ as $PP\left[f(s);s=\lambda_k\right]$, we have, if $p<0$.
\begin{equation*}
       f(z) = \sum_{k=0}^{\infty} \operatorname{PP}(f(z); z = \lambda_k).
\end{equation*}
%\begin{corollary}\label{cor1D}
%As a result
%\begin{equation*}
%   f(s)=\sum_{k=-\infty}^{\infty}\dfrac{a_k}{s-\lambda_k}
%\end{equation*}
%where $a_k$ is the residual of $f(s)$ at $s=\lambda_k$, i. e.
%$$  a_k=\dfrac{(-1)^k}{2T_d}.$$
%\end{corollary}
%\textbf{Proof:}
%\begin{eqnarray*}
%   a_k&=&\dfrac{1}{2j\pi}\oint_{\mathcal{C}_k}f(s)ds=\lim_{s\rightarrow \lambda_k}(s-\lambda_k)f(s)=
%   \lim_{s\rightarrow \lambda_k}\dfrac{(s-\lambda_k)e^{-Tds}}{1-e^{-2T_ds}}\\
%   &=&
%   \lim_{s\rightarrow \lambda_k}\dfrac{e^{-sT_d}-(s-\lambda_k)T_de^{-sT_d}}{2T_de^{-2sT_d}}=\dfrac{e^{\lambda_kT_d}}{2T_d}=\dfrac{e^{\frac{kj\pi}{T_d}T_d}}{2Td}=
%   \dfrac{e^{jk\pi}}{2T_d}=\dfrac{(-1)^k}{2T_d}
%\end{eqnarray*}
%
%$\hfill\Box$

%%%%%%%%%%%%%%%%%%%%%%%%%%%%%%%%


\begin{thebibliography}{99}
\bibitem{Bernstein}
D. S. Bernstein, {\it  Matrix Mathematics - Theory, Facts, and Formulas with Application
to Linear System Theory}, Princeton University Press, Princeton NJ; 2005.

\bibitem{Carvalho:1993}
J. L. Martins de Carvalho, {\it Dynamical Systems and Automatic Control}, Prentice Hall, London; 1999.
\bibitem{SimLiw:2004}
Simone Research Group and Liwacom, {\em  Simone software: equations and methods}, Simone Research Group and LIWACOM, Germany, 2004.

\bibitem{Yueh:2005}
W-C Yueh, "Eigenvalues of Several Tridiagonal Matrices", in  {\it  Applied Mathematics E-Notes},
Vol 5, pp. 72, 2005.

\end{thebibliography}
\end{document}